# A peculiar modular form of weight one

by

Stephen S. Kudla[1]

Michael Rapoport

and

Tonghai Yang[2]

**Introduction.**

In this paper we construct a modular form $\phi$ of weight one, attached to an imaginary quadratic field, which has several curious properties. Let $\boldsymbol{k} = \mathbb{Q}(\sqrt{-q})$ where $q > 3$ is a prime congruent to $3$ modulo $4$. The first curious property concerns the Fourier expansion of $\phi$, which is nonholomorphic and is not a cusp form. For a positive integer $n$, the $n$th Fourier coefficient $a_n(\phi)$ of $\phi$ vanishes unless either $n$ or $n/p$, for some prime $p$, inert in $\boldsymbol{k}$, is a norm of an integral ideal of $\boldsymbol{k}$. If $n$ or $n/p$ is a norm, then $a_n(\phi)$ has a simple formula. The nonholomorphic part of $\phi$ occurs in the negative Fourier coefficients, which are nonvanishing only for integers $-n$ where $n \geq 1$ is a norm and involve the exponential integral. The second curious property concerns the Mellin transform of $\phi$. It has a simple expression in terms of the Dedekind zeta function of $\boldsymbol{k}$ and the difference of the logarithmic derivatives of the Riemann zeta function and the Dirichlet L-series associated to $\boldsymbol{k}$. The third curious property concerns the Fourier coefficients $a_n(\phi)$, for positive $n$, which are connected with the theory of complex multiplication. More precisely, they are the degrees of certain $0$-cycles on the moduli scheme of elliptic curves with complex multiplication by $\mathcal{O}_{\boldsymbol{k}}$. We now describe these results in more detail.

In [**13**], a certain family of (incoherent) Siegel Eisenstein series was introduced which have an odd functional equation and hence have a natural zero at their center of symmetry ($s = 0$). It was suggested that the derivatives of such series at $s = 0$ should have some connection with arithmetical algebraic geometry, and some evidence was provided in the case of genus $2$ and weight $3/2$. In that case, certain of the Fourier coefficients of the central derivative were shown to involve (parts of) the height pairing of Heegner points on Shimura curves. Higher dimensional cases are studied in [**15**] and [**17**].

In the present paper, we consider the simplest possible example of an incoherent

[1] Partially supported by NSF Grant DMS-9622987
[2] Partially supported by NSF Grant DMS-9700777

Typeset by $\mathcal{A}_{\mathcal{M}}\mathcal{S}$-TEX

Eisenstein series and its central derivative. More precisely, let $q > 3$ be a prime congruent to $3$ modulo $4$. Associated to the imaginary quadratic field $\boldsymbol{k} = \mathbb{Q}(\sqrt{-q})$, there is a nonholomorphic Eisenstein series of weight $1$ and character $\chi_q(\cdot) = \left(\frac{-q}{\cdot}\right)$ for the congruence subgroup $\Gamma_0(q)$ of $\Gamma = SL_2(\mathbb{Z})$. In classical language, this series has the following form. For $\tau = u+iv$ in the upper half plane and $s \in \mathbb{C}$ with $\text{Re}(s) > 1$, let

$$(0.1) \qquad E(\tau, s) = v^{s/2} \sum_{\gamma \in \Gamma_\infty \backslash \Gamma} (c\tau + d)^{-1} |c\tau + d|^{-s} \Phi_q(\gamma),$$

where, for $\gamma = \begin{pmatrix} a & b \\ c & d \end{pmatrix} \in \Gamma$,

$$(0.2) \qquad \Phi_q(\gamma) = \begin{cases} \chi_q(a) & \text{if } c \equiv 0 \mod (q), \\ -iq^{-1/2}\chi_q(c) & \text{if } c \text{ is prime to } q. \end{cases}$$

After analytic continuation to the whole $s$ plane, the normalized Eisenstein series

$$(0.3) \qquad E^*(\tau, s) = q^{(s+1)/2} \Lambda(s+1, \chi_q) E(\tau, s)$$

satisfies the functional equation $E^*(\tau, -s) = -E^*(\tau, s)$. Here

$$(0.4) \qquad \Lambda(s, \chi_q) = \pi^{-(s+1)/2} \Gamma\left(\frac{s+1}{2}\right) L(s, \chi_q).$$

Our 'peculiar' modular form of weight $1$ is then

$$(0.5) \qquad \phi(\tau) := -h_{\boldsymbol{k}} \cdot \frac{\partial}{\partial s}\left\{E(\tau, s)\right\}\bigg|_{s=0} = -\frac{\partial}{\partial s}\left\{E^*(\tau, s)\right\}\bigg|_{s=0},$$

so that, up to sign, $\phi(\tau)$ is the leading term of the Laurent expansion of $E^*(\tau, s)$ at the central point $s = 0$. Here $h_{\boldsymbol{k}}$ is the class number of $\boldsymbol{k}$.

**Theorem 1.** *The nonholomorphic modular form $\phi$ has a Fourier expansion:*

$$\phi(\tau) = a_0(\phi, v) + \sum_{n<0} a_n(\phi, v) e(n\tau) + \sum_{n>0} a_n(\phi) e(n\tau),$$

*where the coefficients are as follows. Let*

$$\rho(n) = |\{\mathfrak{a} \subset \mathcal{O}_{\boldsymbol{k}} \mid N(\mathfrak{a}) = n\}|,$$

*so that*

$$\zeta_{\boldsymbol{k}}(s) = \sum_{n=1}^\infty \rho(n) n^{-s},$$

is the Dedekind zeta function of $k$. The constant term of $\phi$ is given by

$$a_0(\phi, v) = -h_k \cdot \left( \log(q) + \log(v) + 2\frac{\Lambda'(1, \chi_q)}{\Lambda(1, \chi_q)} \right).$$

For $n > 0$,

$$a_n(\phi) = 2\log(q) \left(\mathrm{ord}_q(n) + 1\right) \rho(n) + 2\sum_{p \neq q} \log(p) \left(\mathrm{ord}_p(n) + 1\right) \rho(n/p),$$

where the sum runs over primes $p$ which are inert in $k$. For $n < 0$,

$$a_n(\phi, v) = 2\beta_1(4\pi|n|v) \rho(-n).$$

Here $\beta_1$ is essentially the exponential integral:

$$(0.6) \qquad \beta_1(t) = -\operatorname{Ei}(-t) = \int_1^\infty u^{-1} e^{-ut}\, du.$$

Since one of the aims of this paper is to illustrate the constructions of [13] in the simplest possible case, we work adelically, rather than classically, and obtain both the Fourier expansion of $E(\tau, s)$ and its derivative by computing local Whittaker functions and their derivatives. This is done in sections 1–3.

Following the tradition of Hecke, we next compute the Mellin transform of $\phi$, or, more precisely, of $\phi$ with its constant term omitted:

$$(0.7) \qquad \Lambda(s, \phi) := \int_0^\infty \left( \phi(iv) - a_0(\phi, v) \right) v^{s-1}\, dv.$$

This can be done termwise. Let $\Lambda(s) = \pi^{-s/2} \Gamma(\tfrac{s}{2}) \zeta(s)$, and $\Lambda_k(s) = \Lambda(s)\Lambda(s, \chi_q)$.

**Theorem 2.**
$$\Lambda(s, \phi) = \Lambda_k(s) \left[ \log(q) + \frac{\Lambda'(s, \chi)}{\Lambda(s, \chi)} - \frac{\Lambda'(s)}{\Lambda(s)} \right]$$

What seems striking (to us) here is the occurrence of the *difference* of the logarithmic derivatives of $\Lambda(s, \chi)$ and $\Lambda(s)$. The computation is done in section 4.

Next, we show that the Fourier coefficients of $\phi$ are connected with certain 0-cycles on an arithmetic curve. This arithmetic curve is defined as the moduli scheme $M$ of elliptic curves with complex multiplication $\iota$ by the ring of integers $\mathcal{O}_k$ in $k$. For each positive integer $n$, there is a 0-cycle $Z(n)$ on $M$, which is the locus on which the elliptic curves have an additional endomorphism $y$ with $y^2 = -n$, Galois commuting with the action of $\mathcal{O}_k$, i.e, with $y \cdot \iota(a) = \iota(\bar{a}) \cdot y$, for $a \in \mathcal{O}_k$.

The following result is proved in section 5 as an application of the classical theory of complex multiplication and its further development due to Gross and Zagier, [10], [11].

**Theorem 3.** *(i) For $n \in \mathbb{Z}_{>0}$, suppose that neither $n$ nor $n/p$, for any prime $p$ which is inert in $\mathbf{k}$ is a norm of an integral ideal of $\mathcal{O}_{\mathbf{k}}$, i.e., suppose that $\rho(n) = 0$ and $\rho(n/p) = 0$ for all inert primes $p$. Then the cycle $Z(n)$ is empty and $a_n(\phi) = 0$.*

*(ii) If $\rho(n) \neq 0$ ( resp. $\rho(n/p) \neq 0$ ), then the $0$-cycle $Z(n)$ is supported in the fiber at $q$ ( resp. the fiber at $p$ ), and*

$$\deg(Z(n)) = a_n(\phi).$$

We obtain the following Fourier expansion:

$$(0.8) \qquad \phi(\tau) = a_0(\phi, v) + \sum_{n>0} \deg(Z(n)) \, q^n + 2 \sum_{n>0} \rho(n) \, \beta_1(4\pi n v) q^{-n},$$

in which, for a moment, we write, $q^n = e(n\tau) = e^{2\pi i n \tau}$. As pointed out above, this expression provides some evidence in favor of the general proposals of [13]. In section 6, we propose a 'modular' definition of an Arakelov divisor $Z(t, v)$ for $t < 0$ for which

$$\deg(Z(t, v)) = a_t(\phi, v)$$

as well. It would be nice to have a geometric interpretation of the constant term $a_0(\phi, v)$.

Also in section 6, we suggest the modifications which are necessary to extend Theorems 1 and 3 to the most general incoherent Eisenstein series on $SL(2)$.

As far as we know, the modular form $\phi$ does not appear in the classical literature. It resembles, however, the modular form of weight $3/2$ considered by Zagier [24], whose Fourier expansion is

$$(0.9) \qquad \phi_{\text{Zagier}}(\tau) = -\frac{1}{12} + \sum_{n>0} \deg(Z(n)) q^n + v^{-1/2} \sum_{m \in \mathbb{Z}} \beta_{3/2}(4\pi m^2 v) q^{-m^2},$$

where

$$(0.10) \qquad \beta_{3/2}(t) = \frac{1}{16\pi} \int_1^\infty u^{-3/2} e^{-ut} \, du.$$

Here we have momentarily changed notation and have denoted by $Z(n)$ certain 0-cycles on a modular curve (over $\mathbb{C}$), cf., for example, [1].

We conclude with the observation of Zagier that the coefficients of $\phi$ actually occur in [11]. Let $-d$ be a fundamental discriminant with $(d, q) = 1$ and with $d > 4$.

Let

(0.11) $$J(-q,-d) = \prod_{\substack{[\tau_1],[\tau_2] \\ disc(\tau_1) = -q \\ disc(\tau_2) = -d}} \bigl(j(\tau_1) - j(\tau_2)\bigr)$$

be the integer of (1.2) of [**11**]. Then Theorem 1.3 and (1.4) of that paper can be written as

(0.12) $$2\log(J(-d,-q)) = \sum_{\substack{n \\ n^2 < dq \\ n^2 \equiv dq(4)}} a_m(\phi),$$

where, in the sum, $m = (dq - n^2)/4$.

This paper has its origin in a question raised by Dick Gross at the Durham conference in the summer of 1996 when we presented results on the higher dimensional cases. We thank him for his insistence that we consider this 'simplest' case.

**Notation.** We will use the following notation for elements of $SL(2)$:

$$m(a) = \begin{pmatrix} a & \\ & a^{-1} \end{pmatrix}, \ n(b) = \begin{pmatrix} 1 & b \\ & 1 \end{pmatrix}, \ w = \begin{pmatrix} & 1 \\ -1 & \end{pmatrix}, \ k_\theta = \begin{pmatrix} \cos(\theta) & \sin(\theta) \\ -\sin(\theta) & \cos(\theta) \end{pmatrix},$$

for $\theta \in \mathbb{R}$.

Let $\psi : \mathbb{Q}_\mathbb{A} \to \mathbb{C}^1$ be the standard nontrivial additive character, as defined in Tate's thesis.

**§1. Incoherent Eisenstein series for $SL(2)$.**

The simplest examples of incoherent Eisenstein series, as considered in [**13**], occurs for the group $SL(2)$. In this section, we sketch their construction.

Let $\chi$ be a nontrivial quadratic character of $\mathbb{Q}^\times \backslash \mathbb{Q}_\mathbb{A}^\times$, and let $\boldsymbol{k} = \mathbb{Q}(\sqrt{D})$ be the associated quadratic extension. We assume that $D < 0$ is the discriminant of the extension $\boldsymbol{k}/\mathbb{Q}$, and write $\mathcal{O}_{\boldsymbol{k}}$ for the ring of integers of $\boldsymbol{k}$.

Let $V$ be a two dimensional vector space over $\mathbb{Q}$ with quadratic form $x \mapsto Q(x)$. Let $\chi_V(x) = (x, -\det V)$ be the quadratic character associated to $V$, and assume that $\chi_V = \chi$. Then $V \simeq \boldsymbol{k}$ with quadratic form $Q(x) = \kappa \cdot N(x)$, $\kappa \in \mathbb{Q}^\times$, given by a multiple of the norm form. Up to isomorphism over $\mathbb{Q}$, the space $V$ is determined by its collection of localizations $\{V_p\}$, for $p \leq \infty$. These localizations

$V_p$ are determined by their Hasse invariants. They have the following form, up to isomorphism:

(i) If $p = \infty$, $V_\infty$ has signature $(2, 0)$ or $(0, 2)$, and these spaces have Hasse invariants $\epsilon_\infty(V_\infty) = 1$ and $\epsilon_\infty(V_\infty) = -1$, respectively.

(ii) If $p$ is inert or ramified in $k/\mathbb{Q}$, then $V_p = k_p$ with $Q(x) = N(x)$ or $Q(x) = \kappa_p \cdot N(x)$, where $\chi_p(\kappa_p) = -1$, i.e., with $\kappa_p \notin N k_p^\times$. These spaces have Hasse invariants $\epsilon_p(V_p) = 1$ and $\epsilon_p(V_p) = -1$, respectively.

(iii) If $p$ splits in $k$, $V_p = V_{1,1}$, the split binary quadratic space. In this case, the Hasse invariant is $\epsilon_p(V_p) = 1$.

If $V = k$ with $Q(x) = \kappa \cdot N(x)$, then $\epsilon_p(V_p) = \chi_p(\kappa) = (\kappa, D)_p$. The Hasse invariants satisfy the product formula $\prod_{p \leq \infty} \epsilon_p(V_p) = 1$.

An *incoherent collection* $\mathcal{C} = \{\mathcal{C}_p\}$ is a collection of such local quadratic spaces with the following property. For some (and hence for any) global binary quadratic space $V$ with $\chi_V = \chi$, there is a finite set $S$ of places, with (i) $|S|$ is odd, (ii) for any finite place $p \in S$, $p$ does not split in $k$, and (iii) for all $p \leq \infty$,

$$(1.1) \qquad \mathcal{C}_p = \begin{cases} (V_p, \kappa_p Q) & \text{if } p \in S, \\ (V_p, Q) & \text{otherwise.} \end{cases}$$

Here $\kappa_p$ is as in (ii) above and $\kappa_\infty = -1$. By the definition $\prod_{p \leq \infty} \epsilon_p(\mathcal{C}_p) = -1$, and there is no global binary quadratic space with these localizations.

Let $G = SL(2)$ over $\mathbb{Q}$, and let $B = TN$ be the upper triangular Borel subgroup. Associated to quadratic character $\chi$ is the global induced representation

$$I(s, \chi) = \text{Ind}_{B(\mathbb{A})}^{G(\mathbb{A})}(\chi | \ |^s),$$

of $G(\mathbb{A})$. Here $\Phi(s) \in I(s, \chi)$ satisfies

$$(1.2) \qquad \Phi(n(b)m(a)g, s) = \chi(a)|a|^{s+1}\Phi(g, s),$$

for $b \in \mathbb{A}$ and $a \in \mathbb{A}^\times$. At $s = 0$, we have a decomposition into irreducible representations of $G(\mathbb{A})$:

$$(1.3) \qquad I(0, \chi) = \left( \oplus_V \Pi(V) \right) \oplus \left( \oplus_\mathcal{C} \Pi(\mathcal{C}) \right).$$

Here $V$ runs over binary quadratic spaces with $\chi_V = \chi$, and $\mathcal{C}$ runs over incoherent collections with $\chi_\mathcal{C} = \chi$. Let $S(V(\mathbb{A}))$ be the Schwartz space of $V(\mathbb{A})$. Recall that the irreducible subspace $\Pi(V)$ is the image of the map $\lambda_V : S(V(\mathbb{A})) \to I(0, \chi)$ defined by $\lambda_V(\varphi)(g) = (\omega(g)\varphi)(0)$, where $\omega = \omega_\psi$ denotes the action of $G(\mathbb{A})$ on $S(V(\mathbb{A}))$ via the global Weil representation determined by the fixed additive character $\psi$. Analogously, let $S(\mathcal{C}_{\mathbb{A}}) = \otimes_{p \leq \infty} S(\mathcal{C}_p)$ be the Schwartz space of $\mathcal{C}$ (restricted tensor product of the spaces $S(\mathcal{C}_p)$). The group $G(\mathbb{A})$ acts on this space by the restricted product of the local Weil representations, and $\Pi(\mathcal{C})$ is the image of the equivariant map $\lambda_{\mathcal{C}} : S(\mathcal{C}_{\mathbb{A}}) \to I(0, \chi)$ defined by

$$(1.4) \qquad \lambda_{\mathcal{C}}(\otimes_p \varphi_p)(g) = \otimes_p \big((\omega_p(g)\varphi_p)(0)\big),$$

where $\omega_p = \omega_{\psi_p}$ denotes the action of $G(\mathbb{Q}_p)$ on $S(\mathcal{C}_p)$ via the local Weil representation determined by the local component $\psi_p$ of $\psi$.

For $g \in G(\mathbb{A})$, write $g = n(b)m(a)k$, where $k \in K_\infty K$, $K_\infty = SO(2)$, and $K = SL_2(\hat{\mathbb{Z}})$. Here $a \in \mathbb{A}^\times$ is not uniquely determined, but the absolute value $|a(g)| := |a|_{\mathbb{A}}$ is well defined. For $\Phi \in I(0, \chi)$, let

$$(1.5) \qquad \Phi(g, s) = \Phi(g, 0)|a(g)|^s$$

be its standard extension. For such a standard section $\Phi(s)$, the Eisenstein series

$$(1.6) \qquad E(g, s, \Phi) = \sum_{\gamma \in B(\mathbb{Q}) \backslash G(\mathbb{Q})} \Phi(\gamma g, s)$$

converges absolutely for $\mathrm{Re}(s) > 1$ and has an entire analytic continuation.

If $\Phi = \Phi(0) \in \Pi(\mathcal{C})$ for some incoherent collection $\mathcal{C}$, the resulting incoherent Eisenstein series $E(g, s, \Phi)$ vanishes at $s = 0$, [**13**], Theorem 2.2, and [**14**], Theorem 3.1. In the sections which follow, we will compute the central derivative $E'(g, 0, \Phi)$ of such a series in the simplest possible case.

## §2. The case of prime discriminant.

Fix a prime $q > 3$ with $q \equiv 3 \mod 4$, and let $\boldsymbol{k} = \mathbb{Q}(\sqrt{-q})$. Let $R = \mathcal{O}_{\boldsymbol{k}}$ be the ring of integers of $\boldsymbol{k}$. Let $\chi$ be the character of $\mathbb{Q}_{\mathbb{A}}^\times$ associated to $\boldsymbol{k}$, so that, $\chi(x) = (x, -q)_{\mathbb{A}}$, where $(\,,\,)_{\mathbb{A}}$ is the global Hilbert symbol. For a prime $p$, let $\chi_p(x) = (x, -q)_p$, for the local Hilbert symbol $(\,,\,)_p$ at $p$.

Let $V = \boldsymbol{k}$, viewed as a $\mathbb{Q}$-vector space with quadratic form $Q(x) = -N(x)$, and let $\mathcal{C}$ be the incoherent collection defined by $\mathcal{C}_p = V_p$, for all finite primes $p$, and $\mathrm{sig}(\mathcal{C}_\infty) = (2, 0)$. We identify $\mathcal{C}_\infty$ with $\boldsymbol{k}_\infty$ with quadratic form $N(x)$.

We will consider the incoherent Eisenstein series associated to the standard factorizable section $\Phi(s) \in I(s, \chi)$ with $\Phi(0) = \lambda(\varphi)$ for $\varphi = \otimes_p \varphi_p \in S(\mathcal{C}_{\mathbb{A}})$ and $\lambda = \lambda_{\mathcal{C}} : S(\mathcal{C}_{\mathbb{A}}) \to I(0, \chi)$, as above. We choose $\varphi = \otimes_{p \leq \infty} \varphi_p \in S(\mathcal{C}_{\mathbb{A}})$ as follows. For $p < \infty$, let $\varphi_p$ be the characteristic function of $R_p = R \otimes_{\mathbb{Z}} \mathbb{Z}_p$, and let $\varphi_{\infty}(x) = \exp(-\pi N(x))$.

In the range of absolute convergence, $\text{Re}(s) > 1$, the incoherent Eisenstein series determined by $\Phi(s)$ has a Fourier expansion

$$(2.1) \qquad E(g, s, \Phi) = \sum_{t \in \mathbb{Q}} E_t(g, s, \Phi),$$

with

$$(2.2) \qquad E_t(g, s, \Phi) = \int_{\mathbb{Q} \backslash \mathbb{A}} E(n(b)g, s, \Phi) \, \psi(-tb) \, db.$$

For $t \neq 0$,

$$(2.3) \qquad E_t(g, s, \Phi) = \prod_{p \leq \infty} W_{t,p}(g_p, s, \Phi_p),$$

where

$$(2.4) \qquad W_{t,p}(g_p, s, \Phi_p) = \int_{\mathbb{Q}_p} \Phi_p(wn(b)g_p, s) \, \psi_p(-tb) \, db$$

is the local Whittaker integral. Here $db$ is the self dual measure with respect to $\psi_p$. Similarly,

$$(2.5) \qquad E_0(g, s, \Phi) = \Phi(g, s) + M(s)\Phi(g),$$

where the global intertwining operator has a factorization $M(s)\Phi = \otimes_p M_p(s)\Phi_p$, and

$$(2.6) \qquad M_p(s)\Phi(g) = W_{0,p}(g, s, \Phi).$$

To obtain an explicit formula for our particular choice of $\varphi \in S(\mathcal{C}_{\mathbb{A}})$ it is necessary, first, to determine the section $\Phi(s) = \otimes_p \Phi_p(s)$ defined by $\varphi$ and then to compute the integrals $W_{t,p}(g_p, s, \Phi_p)$. For $p \neq q$, these are completely standard matters, and we will spare the reader the details of the computations and just summarize the results. Also, in the end, we will assume that, for $\tau = u + iv$ in the upper half plane, $g = g_{\tau} = n(u)m(v^{1/2}) \in G(\mathbb{R})$, so that we may as well assume that $g_p = 1$ for all finite places. Since our section $\Phi(s)$ has been fixed, we will omit it from the notation and simply write:

$$(2.7) \qquad W^*_{t,p}(s) = L_p(s+1, \chi) \, W_{t,p}(e, s, \Phi_p),$$

for $p$ a finite prime, and

$$(2.8) \qquad W^*_{t,\infty}(\tau, s) = L_{\infty}(s+1, \chi) \, W_{t,\infty}(g_{\tau}, s, \Phi_{\infty}).$$

**Lemma 2.1.** *(i) For $p \neq q$, $\infty$, the function $\varphi_p$ is invariant under $K_p = SL_2(\mathbb{Z}_p)$, and so $\Phi_p(s)$ is the unique $K_p$-invariant vector in $I_p(s,\chi)$ with $\Phi_p(e,s) = 1$.*
*(ii) For $p = \infty$, and for $k_\theta \in SO(2) = K_\infty$, $\omega(k_\theta)\varphi_\infty = e^{i\theta} \cdot \varphi_\infty$. Thus, $\Phi_\infty(s)$ is the unique weight $1$ eigenvector for $K_\infty$ in $I_\infty(s, \chi_\infty)$ with $\Phi_\infty(e,s) = 1$.*

For $p = q$, let $J_q \subset K_q = SL_2(\mathbb{Z}_q)$ be the Iwahori subgroup. The ramified character $\chi_q$, $\chi_q(t) = (t, -q)_q$ defines a character of $J_q$ by

$$\chi_q(\begin{pmatrix} a & b \\ qc & d \end{pmatrix}) = \chi_q(a). \tag{2.9}$$

Recall that $K_q = J_q \cup (J_q \cap N) w J_q$. Thus, the subspace of $I_q(s, \chi)$ consisting of $\chi_q$-eigenvectors of $J_q$ is two dimensional, and is spanned by the 'cell' functions $\Phi_q^0(s)$ and $\Phi_q^1(s)$, determined by

$$\Phi_q^i(w_j, s) = \delta_{ij}, \qquad \text{where } w_0 = 1 \text{ and } w_1 = w. \tag{2.10}$$

By a direct calculation, given in section 3, we find:

**Lemma 2.2.** *For $k = \begin{pmatrix} a & b \\ qc & d \end{pmatrix} \in J_q$, $\omega(k)\varphi_q = \chi_q(a)\varphi_q$, so that $\Phi_q(s)$ lies in the $2$ dimensional subspace of $I_q(s,\chi)$ consisting of $\chi_q$ eigenvectors of $J_q$. Explicitly $\Phi_q(s) = \Phi_q^0(s) + c_q \cdot \Phi_q^1(s)$, where $c_q = \sqrt{-1}\, q^{-1/2}$.*

It follows that our section $\Phi(s)$ is a $\chi$-eigenvector under the action of the group $J = J_q \prod_{p \neq q} K_p$, and hence that the Eisenstein series $E(g, s, \Phi)$ is right $\chi$-equivariant under $J$.

**Lemma 2.3.** *If $g \in G(\mathbb{R})$, then $E_t(g, s, \Phi) = 0$ for $t \notin \mathbb{Z}$.*

*Proof.* For $t \in \mathbb{Q} - \mathbb{Z}$, take $b \in \hat{\mathbb{Z}}$ so that $\psi(tb) \neq 1$. Then, since $n(b) \in J$ with $\chi(n(b)) = 1$, we have

$$E_t(g, s, \Phi) = E_t(gn(b), s, \Phi) = E_t(n(b)g, s, \Phi) = \psi(tb) \cdot E_t(g, s, \Phi), \tag{2.11}$$

and hence $E_t(g, s, \Phi) = 0$. □

From now on we will assume that $t$ is integral.

Next, we compute the normalized local Whittaker functions $W_{t,p}^*(s)$, $W_{t,\infty}^*(\tau, s)$, and their derivatives.

**Lemma 2.4.** *For a finite prime $p \neq q$, let $X = p^{-s}$. Then*

$$W_{t,p}^*(s) = \sum_{r=0}^{\operatorname{ord}_p(t)} (\chi_p(p)X)^r.$$

*Here $\chi_p(p) = (p, -q)_p$ is $1$ if $p$ is split in $\mathbf{k}$ and $-1$ if $p$ is inert. If $t = 0$, so that $\operatorname{ord}_p(t) = \infty$,*

$$W_{0,p}^*(s) = (1 - \chi(p)X)^{-1},$$

*and hence*

$$M_p^*(s)\Phi_p^0 = L_p(s, \chi)\, \Phi_p^0(-s).$$

Here

(2.12) $$M_p^*(s) = L_p(s+1, \chi)\, M_p(s).$$

Let

(2.13) $$\rho_p(t) = \rho(p^{\operatorname{ord}_p(t)}),$$

so that, for $t > 0$,

(2.14) $$\rho(t) = \prod_p \rho_p(t).$$

Note that $\rho(t) = 0$ for $t < 0$.

The following facts, which follow immediately from Lemma 2.4, will be useful below.

**Lemma 2.5.** *At $s = 0$,*

$$W_{t,p}^*(0) = \rho_p(t).$$

*If $\chi_p(p) = -1$ and $\operatorname{ord}_p(t)$ is odd, i.e., if $\rho_p(t) = 0$, then*

$$W_{t,p}^{*,\prime}(0) = \log(p)\, \frac{1}{2}(\operatorname{ord}_p(t) + 1)\, \rho_p(t/p).$$

*Note that, $\rho_p(t/p) = 1$ here.*

Next consider the archimedean factor. In section 3, we indicate the proof of the following result.

**Proposition 2.6.** *For $\tau = u + iv$ in the upper half plane, and $g_\tau = n(u)m(v^{1/2})$, as above,*

$$W^*_{t,\infty}(\tau, s) = v^{(1-s)/2} e(tu) \frac{2i\pi^{s/2} e^{2\pi tv}}{\Gamma(s/2)} \int_{\substack{u>0 \\ u>2tv}} e^{-2\pi u} u^{s/2} (u - 2tv)^{s/2-1} du.$$

*This formula can be found in Siegel [21], eq.(14), p. 88. For $t \neq 0$, it extends to an entire function of $s$.*

(i) *If $t > 0$, then $W^*_{t,\infty}(\tau, 0) = 2i \cdot v^{1/2} e(t\tau)$.*

(ii) *If $t = 0$, then*

$$W^*_{0,\infty}(\tau, s) = \left(M^*_\infty(s)\Phi_\infty\right)(g_\tau) = i \cdot v^{(1-s)/2} \cdot \pi^{-(s+1)/2}\Gamma(\frac{s+1}{2}),$$

*and $W^*_{0,\infty}(\tau, 0) = iv^{1/2}$.*

(iii) *If $t < 0$, then $W^*_{t,\infty}(\tau, 0) = 0$, and*

$$W^{*,\prime}_{t,\infty}(g_\tau, 0, \Phi^1_\infty) = i\, v^{1/2}\, e(t\tau) \cdot \beta_1(4\pi|t|v),$$

*where $\beta_1$ is given by (0.6).*

Finally, in the case $p = q$, and noting that $L_q(s+1, \chi) = 1$, a direct calculation given in section 3 below yields:

**Proposition 2.7.** *For $p = q$, let $c_q = \sqrt{-1}q^{-1/2}$, as in Lemma 2.2 above. Then*

$$W^*_{t,q}(s) = \left(1 - \chi_q(t) q^{-s(\mathrm{ord}_q(t)+1)}\right) \cdot c_q.$$

*If $\chi_q(t) = -1$, then $W^*_{t,q}(0) = 2c_q \cdot \rho_q(t)$. If $\chi_q(t) = 1$, then $W^*_{t,q}(0) = 0$, and*

$$W^{*,\prime}_{t,q}(0) = c_q \log(q) \left(\mathrm{ord}_q(t) + 1\right) \rho_q(t).$$

*Note that the factor $\rho_q(t) = 1$.*

Allowing $\mathrm{ord}_q(t)$ to go to infinity in $W^*_{t,q}(s)$, gives $W^*_{0,q}(s) = c_q$ and hence

(2.15) $$M^*_q(s)\Phi_q = c_q \cdot \Phi_q(-s).$$

We can now assemble these facts. For $t \neq 0$, and for $\mathrm{Re}(s) > 1$, we have

(2.16) $$E^*_t(g_\tau, s, \Phi) = q^{(s+1)/2} \cdot W^*_{t,\infty}(\tau, s) \cdot W^*_{t,q}(s) \cdot \prod_{p \neq q} W^*_{t,p}(s),$$

Since the factors in the product for $p$ with $\operatorname{ord}_p(t) = 0$ are equal to $1$, the expressions given in Lemmas 2.4, 2.6 and 2.7 provide the entire analytic continuation of the right side of (2.16). At $s = 0$ the factors on the right side vanish as follows:

(i) $W^*_{t,\infty}(\tau, 0) = 0$, if $t < 0$, i.e., if $\chi_\infty(t) = -1$,
(ii) $W^*_{t,q}(0) = 0$ if $\chi_q(t) = 1$, and
(iii) $W_{t,p}(0) = 0$ if $\chi_p(t) = -1$ for $p \neq q$, $\infty$, .

Note that case (iii) occurs only when $p$ is inert in $\boldsymbol{k}$ and $\operatorname{ord}_p(t)$ is odd.

Since

$$(2.17) \qquad \chi_\infty(t)\chi_q(t)\prod_{p\neq q}\chi_p(t) = \chi(t) = 1,$$

an odd number of factors in (2.16) are forced to vanish at $s = 0$, via (i)–(iii). Precisely one factor will vanish when one of the quantities $\rho(-t)$, $\rho(t)$ or $\rho(t/p)$ is nonzero. Specifically, this factor will be

(i) $W^*_{t,\infty}(\tau, 0)$ when $\rho(-t) \neq 0$, i.e., when $\chi_\infty(t) = -1$, $\chi_q(t) = -1$ and $\chi_p(t) = 1$ for $p \neq q$, $\infty$;
(ii) $W^*_{t,q}(0)$ when $\rho(t) \neq 0$, i.e., when $\chi_\infty(t) = 1$, $\chi_q(t) = 1$ and $\chi_p(t) = 1$ for $p \neq q$, $\infty$; and
(iii) $W_{t,p}(0)$ when $\rho(t/p) \neq 0$, i.e., when $\chi_\infty(t) = 1$, $\chi_q(t) = -1$, $\chi_p(t) = -1$ and $\chi_\ell(t) = 1$ for $\ell \neq p$, $q$, $\infty$.

These are the only cases which contribute to the derivative at $s = 0$. The corresponding nonzero Fourier coefficients of $E^{*,\prime}(g_\tau, 0, \Phi)$ can then be obtained by combining the values of normalized local Whittaker functions and their derivatives described in the Lemmas above.

**Proposition 2.8.** *(i) If $\rho(-t) \neq 0$, then*

$$E^{*,\prime}_t(g_\tau, 0, \Phi) = -2v^{1/2}\, e(t\tau) \cdot \beta_1(4\pi|t|v) \cdot \rho(-t).$$

*(ii) If $\rho(t) \neq 0$, then*

$$E^{*,\prime}_t(g_\tau, 0, \Phi) = -2v^{1/2}\, e(t\tau) \cdot \log(q)\, \bigl(\operatorname{ord}_q(t) + 1\bigr) \cdot \rho(t).$$

*(iii) If $\rho(t/p) \neq 0$, then*

$$E^{*,\prime}_t(g_\tau, 0, \Phi) = -2v^{1/2}\, e(t\tau) \cdot \log(p)\, \bigl(\operatorname{ord}_p(t) + 1\bigr) \cdot \rho(t/p).$$

Since the classical Eisenstein series of the introduction is

$$(2.18) \qquad E^*(\tau, s) = v^{-1/2} E^*(g_\tau, s, \Phi),$$

we obtain the expressions of Theorem 1 for the nonconstant Fourier coefficients of $\phi$.

Finally, we compute the constant term and its derivative:

$$E_0^*(g_\tau, s, \Phi) = q^{(s+1)/2}\Lambda(s+1, \chi_q) \cdot \Phi(g_\tau, s) + q^{(s+1)/2}M^*(s)\Phi(g_\tau)$$
$$(2.19) \qquad = q^{(s+1)/2}\Lambda(s+1, \chi_q) \cdot v^{(s+1)/2} + q^{(s+1)/2} \cdot i \cdot c_q \cdot \Lambda(s, \chi_q) \cdot v^{(1-s)/2}$$
$$= q^{(s+1)/2}\Lambda(s+1, \chi_q) \cdot v^{(s+1)/2} - q^{(1-s)/2}\Lambda(1-s, \chi_q) \cdot v^{(1-s)/2}.$$

Thus

$$(2.20) \qquad E_0^{*,\prime}(g_\tau, 0, \Phi) = 2\frac{\partial}{\partial s}\left\{q^{(s+1)/2}\Lambda(s+1, \chi_q) \cdot v^{(s+1)/2}\right\}\bigg|_{s=0}$$
$$= v^{1/2} \cdot h_k \cdot \left(\log(q) + \log(v) + 2\frac{\Lambda'(1, \chi_q)}{\Lambda(1, \chi_q)}\right)$$

This gives the constant term of $\phi$, and finishes the proof of Theorem 1.

## §3. Some computations.

In this section, we provide the details for the computations of the values and derivatives of the Whittaker functions used in section 2, in the case $p = q$ or $\infty$.

We begin with the case $p = q$.

*Proof of Lemma 2.2.* We compute the action of $SL_2(\mathbb{Z}_q)$ on $S(V(\mathbb{Q}_q))$ for the Weil representation. Write $J_q = (J_q \cap N_q^-)(J_q \cap M_q)(J_q \cap N_q)$. For $b \in \mathbb{Z}_q$,

$$(3.1) \qquad \omega(n(b))\varphi_q(x) = \psi_q(-bN(x))\varphi_q(x) = \varphi_q(x),$$

since $x \in R_q$ implies that $N(x) \in \mathbb{Z}_q$, and $\psi_q$ is unramified. Similarly, for $a \in \mathbb{Z}_q^\times$,

$$(3.2) \qquad \omega(m(a))\varphi_q(x) = \chi_q(a)|a|\varphi_q(xa) = \chi_q(a)\varphi_q(x).$$

Finally, note that $n_-(qc) = \begin{pmatrix} 1 & \\ qc & 1 \end{pmatrix} = wn(-qc)w^{-1}$. Then

$$(3.3) \qquad \begin{aligned} \omega(w^{-1})\varphi_q(x) &= -i\int_{V_q} \varphi_q(y)\,\psi_q(-tr(xy^\sigma))\,dy \\ &= -i\int_{R_q} \psi_q(-tr(xy^\sigma))\,dy \\ &= -i\,\mathrm{meas}(R_q)\hat{\varphi}_q(x), \end{aligned}$$

where $\hat{\varphi}_q$ is the characteristic function of the dual lattice $\hat{R}_q = \frac{1}{\sqrt{-q}} R_q$. Then

(3.4) $\quad \omega(n(-qc))\omega(w^{-1})\varphi_q(x) = -i \operatorname{meas}(R_q)\,\psi_q(qcN(x))\,\hat{\varphi}_q(x),$

and so

(3.5)
$$\begin{aligned}
\omega(n_-(qc))\varphi_q(x) &= \operatorname{meas}(R_q) \int_{V_q} \psi_q(qcN(y))\,\psi_q(tr(xy^\sigma))\,\hat{\varphi}_q(y)\,dy \\
&= \operatorname{meas}(R_q) \int_{\hat{R}_q} \psi_q(qcN(y))\,\psi_q(tr(xy^\sigma))\,dy \\
&= \operatorname{meas}(R_q) \int_{\hat{R}_q} \psi_q(tr(xy^\sigma))\,dy \\
&= \operatorname{meas}(R_q)\operatorname{meas}(\hat{R}_q)\,\varphi_q(x).
\end{aligned}$$

Here $\operatorname{meas}(R_q)\operatorname{meas}(\hat{R}_q) = 1$, since the measure $dy$ is self dual with respect to $\psi_q$. The second statement is clear, since

(3.6) $\quad \lambda_q(\varphi_q) = \lambda_q(\varphi_q)(w_0)\Phi_q^0(0) + \lambda_q(\varphi_q)(w_1)\Phi_q^1(0),$

so that

(3.7)
$$\begin{aligned}
c_q &= \lambda_q(\varphi_q)(w) \\
&= \omega(w)\varphi_q(0) \\
&= i\int_{V_q} \varphi_q(y)\,dy \\
&= i\operatorname{meas}(R_q) \\
&= i\,q^{-1/2}.
\end{aligned}$$

$\square$

*Proof of Proposition 2.7.* Recall that $W_{t,q}^*(s) = W_{t,q}(e, s, \Phi_q)$, where, for $\operatorname{Re}(s) > 1$, the Whittaker integral $W_{t,q}(g, s, \Phi_q)$ is as in (2.4). We will actually determine this function completely. It is easy to check that

(3.8) $\quad W_{t,q}(n(b)m(a)g, s, \Phi) = \chi_q(a)|a|^{-s+1}\,\psi_q(tb)\,W_{a^2 t,q}(g, s, \Phi),$

so that, by Lemma 2.2 and the decomposition

(3.9) $\quad G(\mathbb{Q}_q) = N_q M_q w_0 J_q \cup N_q M_q w_1 J_q,$

it will suffice to compute the functions

(3.10) $\quad W_{ij}(s,t) := W_{t,q}(w_i, s, \Phi_q^j) = \int_{\mathbb{Q}_q} \Phi_q^j(wn(b)w_i, s)\,\psi_q(-tb)\,db,$

for $i, j \in \{0, 1\}$ and the cell functions $\Phi_q^j$ as in (2.10). Let

$$(3.11) \qquad W_{ij}(s,t)^{(0)} = \int_{\mathbb{Z}_q} \Phi_q^j(wn(b)w_i, s)\, \psi_q(-tb)\, db,$$

and

$$(3.12) \qquad W_{ij}(s,t)^{(1)} = \int_{\mathbb{Q}_q - \mathbb{Z}_q} \Phi_q^j(wn(b)w_i, s)\, \psi_q(-tb)\, db.$$

For $b \in \mathbb{Z}_q$, $wn(b)w_i \in K_q$, so the integrand in $W_{ij}(s,t)^{(0)}$ is independent of $s$. Then

$$(3.13) \qquad \begin{aligned} W_{0j}(s,t)^{(0)} &= \int_{\mathbb{Z}_q} \Phi_q^j(wn(b), s)\, \psi_q(-tb)\, db \\ &= \Phi_q^j(w, s) \int_{\mathbb{Z}_q} \psi_q(-tb)\, db \\ &= \delta_{1j} \mathrm{char}(\mathbb{Z}_q)(t), \end{aligned}$$

where $\mathrm{char}(\mathbb{Z}_q)$ is the characteristic function of $\mathbb{Z}_q$. Next,

$$(3.14) \qquad W_{1j}(s,t)^{(0)} = \int_{\mathbb{Z}_q} \Phi_q^j(-wn(b)w^{-1}, s)\, \psi_q(-tb)\, db.$$

Since $wn(qx)w^{-1} \in J_q$, the function $b \mapsto \Phi_q^j(-wn(b)w^{-1}, s)$ depends only on $b$ modulo $q\mathbb{Z}_q$, and thus the integral vanishes if $t \notin q^{-1}\mathbb{Z}_q$. If $t \in q^{-1}\mathbb{Z}_q$, then

$$(3.15) \qquad W_{1j}(s,t)^{(0)} = q^{-1} \chi_q(-1) \sum_{b \in \mathbb{Z}/q\mathbb{Z}} \Phi_q^j(n_-(-b))\, \psi_q(-tb).$$

Note that $n_-(-b) \in J_q$ if and only if $b \in q\mathbb{Z}_q$. If $b \in \mathbb{Z}_q^\times$, then write

$$(3.16) \qquad n_-(-b) = n(x)wk = n(x)w \begin{pmatrix} A & B \\ qC & D \end{pmatrix}.$$

This yields

$$(3.17) \qquad (0,1)n_-(-b) = (-b, 1) = (-1, 0)k = (-A, -B),$$

and thus $A = b$ and $B = -1$. Hence

$$(3.18) \qquad \Phi_q^j(n_-(-b)) = \chi_q(b) \cdot \delta_{1j},$$

and

$$W_{1j}(s,t)^{(0)} = q^{-1}\chi_q(-1)\left(\delta_{0j} + \delta_{1j}\sum_{b\in(\mathbb{Z}/q\mathbb{Z})^\times}\chi_q(b)\,\psi_q(-tb)\right)$$

(3.19)
$$= \begin{cases} q^{-1}\chi_q(-1)\delta_{0j} & \text{if } ord_q(t) \geq 0, \\ q^{-1}\chi_q(-1)\Big(\delta_{0j} + \delta_{1j}\,\chi_q(t)\,\mathfrak{g}(\chi_q)\Big) & \text{if } ord_q(t) = -1. \\ 0 & \text{otherwise.} \end{cases}$$

Here

(3.20)
$$\mathfrak{g}(\chi_q) = \sum_{b\in(\mathbb{Z}/q\mathbb{Z})^\times}\chi_q(b)\,\psi(-b/q)$$

is the Gauss sum for $\chi_q$. We note that $\chi_q(q) = (q,-q)_q = 1$ and $\chi_q(b) = \left(\frac{b}{q}\right)$ for $b \in (\mathbb{Z}/q\mathbb{Z})^\times$, so that $\chi_q(-1) = -1$.

Next consider $W_{ij}(s,t)^{(1)}$. Here, for $ord_q(b) < 0$,

(3.21)
$$wn(b) = \begin{pmatrix} & 1 \\ -1 & -b \end{pmatrix} = \begin{pmatrix} 1 & -b^{-1} \\ & 1 \end{pmatrix}\begin{pmatrix} -b^{-1} & \\ & -b \end{pmatrix}\begin{pmatrix} 1 & \\ b^{-1} & 1 \end{pmatrix}.$$

Thus

(3.22)
$$\Phi_q^j(wn(b)w_i, s) = \chi_q(-b)|b|^{-s-1}\delta_{ij},$$

since $n_-(b^{-1})w_i = w_i w_i^{-1} n_-(b^{-1})w_i$ and

(3.23)
$$w_i^{-1}n_-(b^{-1})w_i \in (J_q \cap N_-) \cup (J_q \cap N)$$

for either $i$. Therefore, recalling that $\chi_q(q) = 1$,

$$W_{ij}(s,t)^{(1)} = \chi_q(-1)\delta_{ij}\sum_{r=1}^{\infty}q^{-rs}\int_{\mathbb{Z}_q^\times}\chi_q(b)\,\psi_q(-tb/q^r)\,db$$

(3.24)
$$= \begin{cases} \chi_q(-1)\,\delta_{ij}\,q^{-1-s(ord_q(t)+1)}\,\chi_q(t)\,\mathfrak{g}(\chi_q) & \text{if } ord_q(t) \geq 0, \\ 0 & \text{otherwise.} \end{cases}$$

Here we use the fact that

(3.25)
$$\int_{\mathbb{Z}_q^\times}\chi_q(b)\,\psi_q(-tb/q^r)\,db = \begin{cases} 0 & \text{if } r \neq ord_q(t)+1, \\ q^{-1}\,\chi_q(t)\,\mathfrak{g}(\chi_q) & \text{if } r = ord_q(t)+1. \end{cases}$$

These results can be summarized as follows: If $ord_q(t) \geq 0$, then, writing $r = ord_q(t) + 1$,

(3.26)
$$\Big(W_{ij}(s,t)\Big) = q^{-1}\chi_q(-1)\begin{pmatrix} q^{-sr}\,\chi_q(t)\,\mathfrak{g}(\chi_q) & q\chi_q(-1) \\ 1 & q^{-sr}\,\chi_q(t)\,\mathfrak{g}(\chi_q) \end{pmatrix}.$$

If $\mathrm{ord}_q(t) = -1$, then

(3.27) $$\Bigl(W_{ij}(s,t)\Bigr) = q^{-1}\chi_q(-1)\begin{pmatrix} 0 & 0 \\ 1 & \chi_q(t)\mathfrak{g}(\chi_q) \end{pmatrix}.$$

Otherwise $(W_{ij}(s,t)) = 0$.

Finally, recall that $\Phi_q(s) = \Phi_q^0(s) + c_q\Phi_q^1(s)$, with $c_q = \sqrt{-1}\,q^{-1/2}$, and that $\mathfrak{g}(\chi_q) = \sqrt{-1}\,q^{1/2}$. Thus, if $\mathrm{ord}_q(t) \geq 0$,

(3.28) $$\begin{pmatrix} W_{t,q}(w_0, s, \Phi_q) \\ W_{t,q}(w_1, s, \Phi_q) \end{pmatrix} = \Bigl(W_{ij}(s,t)\Bigr) \cdot \begin{pmatrix} 1 \\ c_q \end{pmatrix}$$
$$= (1 - \chi_q(t)q^{-s(\mathrm{ord}_q(t)+1)})\begin{pmatrix} c_q \\ -q^{-1} \end{pmatrix}.$$

If $\mathrm{ord}_q(t) = -1$, then

(3.29) $$\begin{pmatrix} W_{t,q}(w_0, s, \Phi_q) \\ W_{t,q}(w_1, s, \Phi_q) \end{pmatrix} = (1 - \chi_q(t))\begin{pmatrix} 0 \\ -q^{-1} \end{pmatrix}.$$

This completes the proof of Proposition 2.7, and a little more. $\square$

We next turn to the case $p = \infty$. Of course, these calculations are rather well known, but we include them for convenience.

*Proof of Proposition 2.6.* Recall that, by (2.8),

(3.30) $$W^*_{t,\infty}(\tau, s) = v^{\frac{1}{2}(1-s)}\, e(tu)\, W^*_{tv,\infty}(e, s)$$
$$= v^{\frac{1}{2}(1-s)}\, e(tu)\, L_\infty(s+1, \chi)\, W(s, tv),$$

where

(3.31) $$L_\infty(s+1, \chi) = \pi^{-s/2-1}\Gamma\Bigl(\frac{s}{2}+1\Bigr),$$

as in (0.4), and where we set

(3.32) $$W(s,t) := W_{t,\infty}(e, s, \Phi_\infty).$$

A simple computation of the Iwasawa decomposition of $wn(b)$ yields

(3.33) $$\Phi_\infty^1(wn(b), s) = i \cdot \frac{(1+ib)}{(1+b^2)^{s/2+1}}.$$

Thus

(3.34)
$$W(s,t) = i \cdot \int_{-\infty}^{\infty} (1+ib)^{-s/2} (1-ib)^{-s/2-1} e(-tb) \, db$$
$$= i \frac{\pi^{s/2+1}}{\Gamma(s/2+1)} \int_{-\infty}^{\infty} (1+ib)^{-s/2} \int_{0}^{\infty} e^{-\pi(1-ib)u} u^{s/2} \, du \, e(-tb) \, db.$$

Recall that

(3.35)
$$\int_{0}^{\infty} e(-xv) \, e^{-xz} \, x^{s-1} \, dx = (z + 2\pi iv)^{-s} \Gamma(s),$$

and hence,

(3.36)
$$\int_{-\infty}^{\infty} e(vx) \, (z + 2\pi iv)^{-s} \, dv = \begin{cases} \frac{1}{\Gamma(s)} \cdot e^{-xz} \, x^{s-1} & \text{if } x > 0, \\ 0 & \text{if } x \leq 0. \end{cases}$$

Thus

(3.37)
$$\int_{-\infty}^{\infty} (1+ib)^{-s/2} \, e(b(\tfrac{1}{2}u - t)) \, db$$
$$= \begin{cases} \frac{(2\pi)^{s/2}}{\Gamma(s/2)} \cdot e^{-2\pi(\tfrac{1}{2}u-t)} (\tfrac{1}{2}u - t)^{s/2-1} & \text{if } \tfrac{1}{2}u - t > 0, \\ 0 & \text{if } \tfrac{1}{2}u - t \leq 0. \end{cases}$$

This gives

(3.38)
$$W(s,t) = i \frac{\pi^{s/2+1}}{\Gamma(s/2+1)} \frac{(2\pi)^{s/2}}{\Gamma(s/2)} \int_{\substack{u>0 \\ u>2t}} e^{-\pi u} u^{s/2} e^{-2\pi(\tfrac{1}{2}u-t)} (\tfrac{1}{2}u - t)^{s/2-1} \, du$$
$$= 2i \cdot \frac{\pi^{s+1} e^{2\pi t}}{\Gamma(s/2)\Gamma(s/2+1)} \int_{\substack{u>0 \\ u>2t}} e^{-2\pi u} u^{s/2} (u - 2t)^{s/2-1} \, du.$$

This yields the first statement of Proposition 2.6.

We are interested in the behavior at $s = 0$. If $t > 0$, the substitution of $u + 2t$ for $u$ yields

(3.39)
$$W(s,t) = 2i \cdot \frac{\pi^{s+1} e^{-2\pi t}}{\Gamma(s/2)\Gamma(s/2+1)} \int_{0}^{\infty} e^{-2\pi u} (u + 2t)^{s/2} u^{s/2-1} \, du$$
$$= 2i \cdot \frac{\pi^{s+1} e^{-2\pi t}}{\Gamma(s/2)\Gamma(s/2+1)} \Big( \Gamma(s/2)(2\pi)^{-s/2}(2t)^{s/2} + O(s) \Big).$$

Thus, if $t > 0$, $W(0, t) = 2\pi i \cdot e^{-2\pi t}$.

If $t = 0$, we get

$$(3.40) \qquad W(s, 0) = 2i \cdot \frac{\pi^{s+1}}{\Gamma(s/2)\Gamma(s/2+1)} \int_0^\infty e^{-2\pi u} u^{s-1} \, du$$

$$= 2\pi i \cdot \frac{2^{-s} \Gamma(s)}{\Gamma(s/2)\Gamma(s/2+1)}$$

$$= i\pi^{1/2} \cdot \frac{\Gamma(\frac{s+1}{2})}{\Gamma(\frac{s+2}{2})},$$

via duplication.

Finally, when $t < 0$,

$$(3.41) \qquad W(s, t) = 2i \cdot \frac{\pi^{s+1} e^{2\pi t}}{\Gamma(s/2)\Gamma(s/2+1)} \int_0^\infty e^{-2\pi u} u^{s/2} (u + 2|t|)^{s/2-1} \, du.$$

The integral here converges for $\operatorname{Re}(s) > -2$, so that the $\Gamma(s/2)$ in the denominator yields $W(0, t) = 0$ when $t < 0$. Note that

$$(3.42) \qquad \begin{aligned} &\int_0^\infty e^{-2\pi u} (u + 2|t|)^{-1} \, du \\ &= \int_0^\infty e^{-4\pi |t| u} (2|t|u + 2|t|)^{-1} 2|t| \, du \\ &= \int_0^\infty \frac{e^{-4\pi |t| u}}{u + 1} \, du \\ &= e^{4\pi |t|} \cdot \int_0^\infty \frac{e^{-4\pi |t| (u+1)}}{u + 1} \, du \\ &= -e^{4\pi |t|} \cdot Ei(-4\pi |t|). \end{aligned}$$

Therefore

$$(3.43) \qquad \left. \frac{\partial}{\partial s} \{ W(s, t) \} \right|_{s=0} = -2\pi i \cdot \frac{1}{2} \cdot e^{2\pi t} \cdot e^{4\pi |t|} \cdot Ei(-4\pi |t|)$$

$$= -i\pi \cdot e^{2\pi |t|} \cdot Ei(-4\pi |t|),$$

when $t < 0$. This completes the proof of Proposition 2.6. $\square$

### §4. The Mellin transform.

In this section, we compute the Mellin transform $\Lambda(s, \phi)$.

**Lemma 4.1.** *The Mellin transform (0.7) can be computed termwise.*

*Proof.* It suffice to check the termwise absolute convergence of the expression

$$-\int_0^\infty \sum_{t=1}^\infty \rho(t) \cdot Ei(-4\pi tv)\, e^{2\pi tv}\, v^{s-1}\, dv. \tag{4.1}$$

It will turn out to be more convenient to consider the expression

$$-\sum_{t=1}^\infty \rho(t) \int_0^\infty Ei(-4\pi tv)\, e^{2\pi tv}\, v^{\sigma-1}\, dv, \tag{4.2}$$

with $\sigma > 1$. First, by a change of variable in each integral, we obtain

$$
\begin{aligned}
&-\sum_{t=1}^\infty \rho(t)(2\pi t)^{-\sigma} \int_0^\infty Ei(-2v)\, e^v\, v^{\sigma-1}\, dv \\
&= -\left(\sum_{t=1}^\infty \rho(t)(2\pi t)^{-\sigma}\right) \cdot \left(\int_0^\infty Ei(-2v)\, e^v\, v^{\sigma-1}\, dv\right) \\
&= -(2\pi)^{-\sigma} \zeta_k(\sigma) \cdot \int_0^\infty Ei(-2v)\, e^v\, v^{\sigma-1}\, dv.
\end{aligned} \tag{4.3}
$$

Thus, it suffices to prove that that the integral is convergent when $\sigma > 1$. Note that the expression

$$-Ei(-2v) = \int_1^\infty \frac{e^{-2vu}}{u}\, du \tag{4.4}$$

is positive. If $v > \delta > 0$, and if $0 < \epsilon < 1$, we can write:

$$
\begin{aligned}
-Ei(-2v)\, e^v &< \int_1^\infty \frac{e^{-(1-\epsilon)vu}}{u}\, du\, e^{-\epsilon v} \\
&< \int_1^\infty \frac{e^{-(1-\epsilon)\delta u}}{u}\, du\, e^{-\epsilon v} \\
&= C_{\delta,\epsilon} \cdot e^{-\epsilon v}.
\end{aligned} \tag{4.5}
$$

On the other hand, if $0 < v < \delta$, we have

$$
\begin{aligned}
-Ei(-2v)\, e^v &< e^\delta \int_1^\infty \frac{e^{-2vu}}{u}\, du \\
&< e^\delta \int_1^\infty e^{-2vu}\, du \\
&= e^\delta (2v)^{-1}.
\end{aligned} \tag{4.6}
$$

Returning to the integral,

$$(4.7) \quad -\int_0^\infty Ei(-2v)\, e^v\, v^{\sigma-1}\, dv < \int_0^\delta e^\delta (2v)^{-1}\, v^{\sigma-1}\, dv + \int_\delta^\infty C_{\delta,\epsilon} \cdot e^{-\epsilon v}\, v^{\sigma-1}\, dv.$$

The second integral on the right side converges for all $\sigma$, while the first is convergent for $\text{Re}(s) > 1$, as required. $\square$

Computing (0.7) termwise, the contribution of the holomorphic terms is simply

$$(4.8) \quad 2(2\pi)^{-s}\Gamma(s) \sum_{t=1}^\infty \left[ \log(q)(\text{ord}_q(t)+1)\,\rho(t) + \sum_{p \neq q,\ \infty} \log(p)(\text{ord}_p(t)+1)\rho(t/p) \right] t^{-s},$$

where the sum on $p$ is over primes inert in $k$. If $p$ is such a prime, then $\rho(p^r) = 1$ if $r$ is even and $0$ if $r$ is odd, and we have

$$\sum_{t=1}^\infty (\text{ord}_p(t)+1)\,\rho(t/p)\, t^{-s} = \sum_{\substack{t=1 \\ (t,p)=1}}^\infty \sum_{r=0}^\infty (r+1)\,\rho(p^{r-1}t)\, p^{-rs} t^{-s}$$

$$(4.9) \quad = \left( \sum_{\substack{t=1 \\ (t,p)=1}}^\infty \rho(t)\, t^{-s} \right) \left( \sum_{r=0}^\infty (r+1)\,\rho(p^{r-1})\, p^{-rs} \right)$$

$$= (1 - p^{-2s}) \cdot \zeta_k(s) \cdot \left( \sum_{r=0}^\infty (2r+2)\, p^{-s(2r+1)} \right)$$

$$= \frac{2p^{-s}}{1 - p^{-2s}} \cdot \zeta_k(s).$$

Similarly,

$$(4.10) \quad \sum_{t=1}^\infty (\text{ord}_q(t)+1)\,\rho(t)\, t^{-s} = (1 - q^{-s})^{-1} \cdot \zeta_k(s).$$

Thus, the holomorphic part of the Fourier expansion contributes

$$(4.11) \quad \Lambda_k(s) \left[ \frac{1}{1-q^{-s}} \log(q) + \sum_{p \neq q} \frac{2p^{-s}}{1-p^{-2s}} \log(p) \right]$$

$$= \Lambda_k(s) \left[ \log(q) + \frac{L'(s,\chi)}{L(s,\chi)} - \frac{\zeta'(s)}{\zeta(s)} \right],$$

as one easily checks.

The termwise transform of the nonholomorphic part of the Fourier expansion is

$$2\sum_{t=1}^{\infty} \rho(t) \int_0^{\infty} \beta_1(4\pi t v)\, e^{2\pi t v}\, v^{s-1}\, dv$$

(4.12)
$$= 2(2\pi)^{-s} \sum_{t=1}^{\infty} \rho(t)\, t^{-s} \left( \int_0^{\infty} \beta_1(2v)\, e^v\, v^{s-1}\, dv \right)$$

$$= 2(2\pi)^{-s}\, \zeta_k(s)\, \Gamma(s) \left[ \frac{\Gamma(\frac{s+1}{2})'}{\Gamma(\frac{s+1}{2})} - \frac{\Gamma(\frac{s}{2})'}{\Gamma(\frac{s}{2})} \right].$$

Thus, the contribution of these terms is

(4.13)
$$\Lambda_k(s) \left[ \frac{\Gamma(\frac{s+1}{2})'}{\Gamma(\frac{s+1}{2})} - \frac{\Gamma(\frac{s}{2})'}{\Gamma(\frac{s}{2})} \right].$$

Combining this with the contribution of the holomorphic terms, we obtain the expression of Theorem 2 for $\Lambda(s, \phi)$.

In (4.12) we have used the following identity, noting that $\Gamma(\frac{s}{2})' = \frac{1}{2}\Gamma'(\frac{s}{2})$.

**Lemma 4.2.**

$$\int_0^{\infty} \beta_1(2v)\, e^v\, v^{s-1}\, dv = \frac{1}{2}\Gamma(s) \left[ \frac{\Gamma'(\frac{s+1}{2})}{\Gamma(\frac{s+1}{2})} - \frac{\Gamma'(\frac{s}{2})}{\Gamma(\frac{s}{2})} \right].$$

*Proof.* Recall that, [**18**],

(4.14)
$$\psi(s) := \frac{\Gamma'(s)}{\Gamma(s)} = -\gamma + \int_0^1 \frac{1 - x^{s-1}}{1 - x}\, dx.$$

Therefore

(4.15)
$$\psi(s + 1/2) - \psi(s) = \int_0^1 \frac{x^{s-1} - x^{s-1/2}}{1 - x}\, dx$$
$$= \int_0^1 x^{s-1}\, \frac{1 - x^{1/2}}{1 - x}\, dx$$
$$= \int_0^1 \frac{x^{s-1}}{1 + x^{1/2}}\, dx$$
$$= 2 \int_0^1 \frac{u^{2s-1}}{1 + u}\, du.$$

Thus we have the useful formula:

(4.16)
$$2 \int_0^1 \frac{w^{s-1}}{1 + w}\, dw = \psi(\frac{s+1}{2}) - \psi(\frac{s}{2}).$$

But now

$$\int_0^\infty \beta_1(2v)\, e^v\, v^{s-1}\, dv$$

(4.17)
$$= \int_0^\infty \int_1^\infty \frac{e^{-2uv}}{u}\, du\, e^v\, v^{s-1}\, dv$$

$$= \Gamma(s) \int_1^\infty (2u-1)^{-s}\, u^{-1}\, du$$

$$= \Gamma(s) \int_1^\infty r^{-s}(1+r)^{-1}\, dr$$

$$= \Gamma(s) \int_0^1 \frac{w^{s-1}}{1+w}\, dw.$$

This gives the claimed expression. □

### §5. Algebraic-geometric aspects.

In the beginning of this section we will consider a slightly more general situation than before. We let $k$ be an imaginary quadratic field with ring of integers $\mathcal{O}_k$. We denote by $a \mapsto \bar{a}$ the non trivial automorphism of $k$.

We will consider the following moduli problem $\mathcal{M}$ over $(\mathrm{Sch}/\mathcal{O}_k)$. To $S \in (\mathrm{Sch}/\mathcal{O}_k)$ we associate the category $\mathcal{M}(S)$ of pairs $(E, \iota)$ where

- $E$ is an elliptic curve over $S$ (i.e. an abelian scheme of relative dimension one over $S$)
- $\iota : \mathcal{O}_k \to \mathrm{End}_S(E)$ is a homomorphism such that the induced homomorphism

(5.1) $$\mathrm{Lie}(\iota) : \mathcal{O}_k \to \mathrm{End}_{\mathcal{O}_S}(\mathrm{Lie}\, E) = \mathcal{O}_S$$

coincides with the structure homomorphism.

The morphisms in this category are the isomorphisms. A pair $(E, \iota)$ will be called an *elliptic curve with CM by $\mathcal{O}_k$*. We denote by $M$ the corresponding set valued functor of isomorphism classes of objects of $\mathcal{M}$.

**Proposition 5.1.** *a) The moduli problem $\mathcal{M}$ is representable by an algebraic stack (in the sense of Deligne-Mumford) which is finite and étale over $\mathrm{Spec}\, \mathcal{O}_k$.*

*b) The functor $M$ has a coarse moduli scheme which is also finite and étale over $\mathrm{Spec}\, \mathcal{O}_k$.*

*Proof.* a) Let $\tilde{\mathcal{M}}$ be the algebraic moduli stack of elliptic curves *over* $\mathrm{Spec}\, \mathcal{O}_k$. The forgetful morphism $i : \mathcal{M} \to \tilde{\mathcal{M}}$ is relatively representable by an unramified morphism (rigidity theorem). Since the generic elliptic curve in any characteristic has only trivial (i.e., $\mathbb{Z}$) endomorphisms, $\mathcal{M}$ has finite fibres over $\operatorname{Spec} \mathcal{O}_k$. Now the finiteness of $\mathcal{M}$ over $\operatorname{Spec} \mathcal{O}_k$ follows from the valuative criterion of properness, since an elliptic curve with complex multiplication has potentially good reduction. The étaleness of $\mathcal{M}$ over $\operatorname{Spec} \mathcal{O}_k$ will follow by the Serre-Tate theorem from a study of the $p$-divisible group of an elliptic curve with $CM$ by $\mathcal{O}_k$.

b) That $M$ has a coarse moduli scheme follows from a) (or use level-$n$-structures). Let $\xi$ be a geometric point of $\mathcal{M}$ and $\bar\xi$ the corresponding geometric point of $M$. Let $\hat{\mathcal{O}}_{\mathcal{M},\xi}$ and $\hat{\mathcal{O}}_{M,\bar\xi}$ be the completions of the local rings. If $\xi$ corresponds to $(E,\iota)$, we have

$$\operatorname{Aut}(\xi) = \operatorname{Aut}(E,\iota) = \mathcal{O}_k^\times , \tag{5.2}$$

comp. cases 1 and 2 below. This group acts on $\hat{\mathcal{O}}_{\mathcal{M},\xi}$ and

$$\hat{\mathcal{O}}_{M,\bar\xi} = (\hat{\mathcal{O}}_{\mathcal{M},\xi})^{\operatorname{Aut}(\xi)} . \tag{5.3}$$

However, since $\mathcal{M}$ is étale over $\operatorname{Spec}\mathcal{O}_k$, we have that $\hat{\mathcal{O}}_{\mathcal{M},\xi}$ is equal to the strict completion of $\mathcal{O}_k$ in the prime ideal corresponding to the image of $\xi$ in $\operatorname{Spec}\mathcal{O}_k$. Hence the action of $\operatorname{Aut}(\xi)$ on $\hat{\mathcal{O}}_{\mathcal{M},\xi}$ is trivial. Therefore the canonical morphism

$$\mathcal{M} \to M \tag{5.4}$$

is étale and the last assertion follows. $\square$

**Corollary 5.2.** *Let $\xi$ be a geometric point of $\mathcal{M}$ and let $\bar\xi$ be the corresponding point of $M$. Then*

$$\hat{\mathcal{O}}_{M,\bar\xi} = \hat{\mathcal{O}}_{\mathcal{M},\xi}.$$

Let $p$ be a prime number and $\wp$ a prime ideal of $\mathcal{O}_k$ above $p$. Let $\kappa$ be an algebraically closed extension of the residue field $\kappa(\wp)$ of $\wp$. Let $(E,\iota) \in \mathcal{M}(\kappa)$, and let $X$ be the $p$-divisible group of $E$. Then the $p$-adic completion $\mathcal{O}_p = \mathcal{O}_k \otimes \mathbb{Z}_p$ acts on $X$. We distinguish two cases.

**Case 1:** *$p$ splits in $\mathcal{O}_k$.* Then $p = \wp \cdot \bar\wp$ with $\wp \neq \bar\wp$ and

$$\mathcal{O}_p = \mathcal{O}_\wp \oplus \mathcal{O}_{\bar\wp} = \mathbb{Z}_p \oplus \mathbb{Z}_p . \tag{5.5}$$

By the condition on $\operatorname{Lie}(\iota)$ the corresponding decomposition of $X$ is of the form

$$X = \hat{\mathbb{G}}_m \times \mathbb{Q}_p/\mathbb{Z}_p . \tag{5.6}$$

In this case it is obvious that $(X, \iota)$ deforms uniquely, i.e. the base of the universal deformation is $\operatorname{Spf} W(\kappa)$, where $W(\kappa)$ denotes the ring of Witt vectors of $\kappa$. In this case the elliptic curve $E$ is ordinary and

(5.7) $$\iota : \mathcal{O}_{\boldsymbol{k}} \xrightarrow{\sim} \operatorname{End}(E) \ .$$

**Case 2:** $p$ *not split in* $\mathcal{O}_{\boldsymbol{k}}$. In this case $\mathcal{O}_p = \mathcal{O}_\wp$ and $(X, \iota)$ is a formal $\mathcal{O}_\wp$-module of height 2 over $\kappa$, in the sense of Drinfeld [**7**], comp. also [**10**]. Let $W_\mathcal{O}(\kappa)$ be the ring of relative Witt vectors of $\kappa$, [**2**]. Then again $(X, \iota)$ deforms uniquely, i.e. the base of the universal deformation of $(X, \iota)$ is $\operatorname{Spf} W_\mathcal{O}(\kappa)$ ([**19**], comp. [**10**], 2.1).

In this case the elliptic curve $E$ is supersingular. Let $\mathbb{B}_p$ be the definite quaternion algebra over $\mathbb{Q}$ which ramifies at $p$. Then $\operatorname{End}(E)$ is isomorphic to a maximal order in $\mathbb{B}_p$.

The main theorem of complex multiplication may be summarized as follows (comp. [**20**],[**3**])

**Proposition 5.3.** *Let $\bar{\boldsymbol{k}}$ be the algebraic closure of $\boldsymbol{k}$. The Galois group $\operatorname{Gal}(\bar{\boldsymbol{k}}/\boldsymbol{k})$ acts on $M(\bar{\boldsymbol{k}})$ through its maximal abelian quotient, which we identify via class field theory with $(\boldsymbol{k} \otimes \mathbb{A}_f)^\times / \boldsymbol{k}^\times$. There is a bijection*

$$M(\bar{\boldsymbol{k}}) = \boldsymbol{k}^\times \backslash (\boldsymbol{k} \otimes \mathbb{A}_f)^\times / (\mathcal{O}_{\boldsymbol{k}} \otimes \hat{\mathbb{Z}})^\times,$$

*such that an element of $\operatorname{Gal}(\bar{\boldsymbol{k}}/\boldsymbol{k})^{\mathrm{ab}}$ represented by an idèle $x \in (\boldsymbol{k} \otimes \mathbb{A}_f)^\times$ acts by translation by $x$ on the right side.*

**Corollary 5.4.** *Let $H$ be the Hilbert class field of $\boldsymbol{k}$ and let $\mathcal{O}_H$ be its ring of integers. Then*

$$M \simeq \operatorname{Spec}(\mathcal{O}_H).$$

**Corollary 5.5.** *Fix $p$ and a prime $\wp$ above $p$. Then $\mathcal{M}(\overline{\kappa(\wp)})$ forms a single isogeny class. There is a bijection*

$$M(\overline{\kappa(\wp)}) = \boldsymbol{k}^\times \backslash (\boldsymbol{k} \otimes \mathbb{A}_f)^\times / (\mathcal{O}_{\boldsymbol{k}} \otimes \hat{\mathbb{Z}})^\times$$

*The action by the Frobenius automorphism over $\kappa(\wp)$ on the left corresponds to the translation by the idèle with component a uniformizer at $\wp$ and 1 at all other finite places on the right.*

**Remark 5.6.** The bijection in Corollary 5.5 is given explicitly as follows. Choose a base point $(E_o, \iota_o) \in \mathcal{M}(\overline{\kappa(\wp)})$, and an identification $\operatorname{End}^0(E_o) = \mathbb{B}_p$. Let $\mathcal{O}(E_o) \subset \mathbb{B}_p$ be the maximal order corresponding to $\mathrm{End}(E_o)$. To $g \in (\boldsymbol{k} \otimes \mathbb{A}_f)^\times$ there is associated the point $(E, \iota) \in \mathcal{M}(\overline{\kappa(\wp)})$ whose Tate module $\hat{T}(E) = \prod_\ell T_\ell(E)$ (with $p$-component the covariant Dieudonné module) is the lattice

$$(5.8) \qquad g \cdot \hat{T}(E_o) \subset \hat{T}(E_o) \otimes \mathbb{Q} \ .$$

In particular

$$(5.9) \qquad \mathcal{O}(E) = \mathbb{B}_p \cap g(\mathcal{O}(E_o) \otimes \hat{\mathbb{Z}}) g^{-1} \ ,$$

is a maximal order in $\mathbb{B}_p$.

We now turn to the definition of special cycles on the moduli space/stack $\mathcal{M}$.

**Definition 5.7.** *Let $(E, \iota) \in \mathcal{M}(S)$. A special endomorphism of $(E, \iota)$ is an element $y \in \mathrm{End}(E)$ with*

$$(5.10) \qquad y \circ \iota(a) = \iota(\overline{a}) \circ y, \qquad \textit{for all } a \in \mathcal{O}_{\boldsymbol{k}} \ .$$

Let $V(E, \iota)$ be the set of special endomorphisms of $(E, \iota)$. Then $V(E, \iota)$ is an $\mathcal{O}_{\boldsymbol{k}}$-module.

Let us first consider the case where $S = \mathrm{Spec}\, \kappa$ with $\kappa$ an algebraically closed field. If $\mathrm{char}\, \kappa = 0$, then $\mathrm{End}(E) = \mathcal{O}_{\boldsymbol{k}}$ and hence $V(E, \iota) = (0)$. The same is true if $\mathrm{char}\, \kappa = p$ and $p$ splits in $\mathcal{O}_{\boldsymbol{k}}$, cf. (5.7). Let now $\mathrm{char}\, \kappa = p$ where $p$ does not split. Then $\mathrm{End}(E)$ can be identified with a maximal order $\mathcal{O}(E)$ in the quaternion algebra $\mathbb{B}_p$ introduced earlier and the $\mathcal{O}_{\boldsymbol{k}}$-module $V(E, \iota)$ is of rank one. Let $a \mapsto a'$ be the main involution of $\mathbb{B}_p$. Then $V(E, \iota)$ is stable under this involution. Let $y \in V(E, \iota)$. Since

$$(5.11) \qquad y + y' \in \mathbb{Z} \cap V(E, \iota) = (0) \ ,$$

we obtain $y^2 = -y \cdot y' \in \mathbb{Z}$. We define a $\mathbb{Z}$-valued quadratic form $Q$ on the $\mathbb{Z}$ module $V(E, \iota)$ of rank 2 by

$$(5.12) \qquad Q(y) = -y^2 = y \cdot y'.$$

The definition of this quadratic form extends to an elliptic curve with $CM$ by $\mathcal{O}_{\boldsymbol{k}}$ over any *connected* scheme $S$.

**Definition 5.8.** *For a positive integer $t$, let $\mathcal{Z}(t)$ be the moduli problem over $\mathcal{M}$ which to $S \in (\mathrm{Sch}/\mathcal{O}_{\boldsymbol{k}})$ associates the category of triples $(E, \iota, y)$ where $(E, \iota) \in \mathcal{M}(S)$ and where $y \in V(E, \iota)$ satisfies*

$$(5.13) \qquad Q(j) = -y^2 = t.$$

We denote by

$$(5.14) \qquad \mathrm{pr} : \mathcal{Z}(t) \longrightarrow \mathcal{M}$$

the forgetful morphism. It is obvious that pr is relatively representable, hence $\mathcal{Z}(t)$ is representable by an algebraic stack. By the rigidity theorem the morphism pr is finite and unramified. We also denote by pr the induced morphism between coarse moduli schemes

$$(5.15) \qquad \mathrm{pr} : Z(t) \longrightarrow M \ .$$

We note that for a geometric point $\tilde{\xi} \in \mathcal{Z}(t)$ we have $\mathrm{Aut}(\tilde{\xi}) = \{\pm 1\}$. Indeed, if $\tilde{\xi}$ corresponds to $(E, \iota, y)$ then an automorphism is given by a unit $u \in \mathcal{O}_{\boldsymbol{k}}^{\times}$ with $uy = yu$, i.e. with $u = \overline{u}$.

**Proposition 5.9.** *The fibre of $Z(t)$ in characteristic 0 or in characteristic $p$ where $p$ splits in $\boldsymbol{k}$ is empty. In particular, $Z(t)$ is an artinian scheme. For $p$ not split in $\boldsymbol{k}$, let $\wp$ be the prime ideal of $\mathcal{O}_{\boldsymbol{k}}$ dividing $p$. The fibre of pr over a geometric point $(E, \iota) \in M(\overline{\kappa(\wp)})$ is empty if $t$ is not represented by the quadratic form on $V(E, \iota)$.*

Let $p$ be inert or ramified in $\boldsymbol{k}$. We choose a base point $(E_o, \iota_o) \in \mathcal{M}(\overline{\kappa(\wp)})$ and denote by $V^{(p)}$ the quadratic space $V(E_o, \iota_o) \otimes_{\mathbb{Z}} \mathbb{Q}$ over $\mathbb{Q}$. Then

$$(5.16) \qquad \mathrm{End}^0(E_o, \iota_o) = \mathbb{B}_p = \boldsymbol{k} \oplus V^{(p)},$$

where the decomposition is orthogonal with respect to the quadratic form given by the reduced norm on $\mathbb{B}_p$. Choosing a set of coset representative $g_i$, $i = 1, \ldots, h_{\boldsymbol{k}}$ for the cosets on the right side of Corollary 5.5 and applying the construction of Remark 5.6, we obtain $\mathcal{O}_{\boldsymbol{k}}$-lattices

$$(5.17) \qquad V(E_i, \iota_i) = V^{(p)} \cap \mathcal{O}(E_i)$$

in $V^{(p)}$.

**Corollary 5.10.** *Let $p$ be inert or ramified in $\mathbf{k}$ and $\wp \subset \mathcal{O}_{\mathbf{k}}$ over $p$. Then*

$$|Z(t)(\overline{\kappa(\wp)})| = \sum_i |\{x \in V(E_i, \iota_i) \mid Q(x) = t\}|.$$

We now quote the result of Gross [**10**], Prop. 3.3 and 4.3. Let $p$ be an inert or ramified prime and $\kappa$ an algebraically closed extension of $\kappa(\wp)$. If $\xi \in \mathcal{M}(\kappa)$, there is an isomorphism

(5.18) $$\hat{\mathcal{O}}_{\mathcal{M},\xi} = W_{\mathcal{O}}(\kappa) ,$$

cf. case 2 above. Here again $W_{\mathcal{O}}(\kappa) = W_{\mathcal{O}_\wp}(\kappa)$ is the ring of relative Witt vectors. We let $\pi$ denote a uniformizer of $\mathcal{O}_\wp$ and hence of $W_{\mathcal{O}}(\kappa)$. Let $\tilde{\xi} \in \mathcal{Z}(t)(\kappa)$ with $\mathrm{pr}(\tilde{\xi}) = \xi$, and $\overline{\tilde{\xi}} \in Z(t)(\kappa)$ its image in $Z(t)$. Here pr is as in (5.14).

**Theorem 5.11.** (Gross): *We have*

$$\hat{\mathcal{O}}_{Z(t),\overline{\tilde{\xi}}} = \hat{\mathcal{O}}_{\mathcal{Z}(t),\tilde{\xi}} = W_{\mathcal{O}}(\kappa)/(\pi^\nu)$$

*where*

$$\nu = \nu_p(t) = \frac{\mathrm{ord}_p(t) + d_p - 1}{f_p} + 1,$$

$$d_p = \mathrm{ord}_p(\mathrm{disc}_{\mathbf{k}/\mathbb{Q}}),$$

*and*

$$f_p = [\kappa(\wp) : \mathbb{F}_p].$$

*In particular, if $p$ is unramified in $\mathbf{k}$, we have $f_p = 2$, and $d_p = 0$ and $\mathrm{ord}_p(t) \equiv 1 (\mathrm{mod}\, 2)$, and*

$$\nu_p(t) = \frac{\mathrm{ord}_p(t) + 1}{2} .$$

*If $p$ is ramified in $\mathbf{k}$, then $f_p = 1$ and $d_p \geq 1$ and, if $p \neq 2$,*

$$\nu_p(t) = \mathrm{ord}_p(t) + 1,$$

*If $p = 2$ is ramified in $\mathbf{k}$, write $\mathbf{k} = \mathbb{Q}(\sqrt{-q})$ with $q$ square free. Then*

$$\nu_p(t) = \begin{cases} \mathrm{ord}_p(t) + 1 & \text{if } q \equiv 3 \mod 4, \\ \mathrm{ord}_p(t) + 2 & \text{if } q \equiv 1 \mod 4, \\ \mathrm{ord}_p(t) + 3 & \text{if } q \equiv 2 \mod 4. \end{cases}$$

Note that the length of $\hat{\mathcal{O}}_{\mathcal{Z}(t),\tilde{\xi}}$ only depends on $\mathrm{ord}_p(t)$. Combining Corollary 5.10 and Theorem 5.11, we therefore obtain an expression for the *degree* of $Z(t)$. Here we define the degree as

(5.19) $$\deg(Z(t)) = \sum_{\xi \in Z(t)} \log\left(\#(\mathcal{O}_{Z(t),\xi})\right),$$

[**22**], [**9**].

**Theorem 5.12.**

$$\deg(Z(t)) = \sum_p f_p \log(p)\, \nu_p(t) \cdot \sum_i |\{x \in V(E_i, \iota_i) \mid Q(x) = t\}|.$$

*Here $p$ ranges over the primes which are not split in $\mathbf{k}$.*

For a fixed $p$ which is not split in $\mathbf{k}$, we now describe the lattices $V(E_i, \iota_i)$ of (5.17) more explicitly. The basic technique is from [**11**], [**6**].

Let $\Delta$ be the discriminant of the field $\mathbf{k}$. Fix an auxillary prime $p_0 \nmid 2p\Delta$ as follows. If $p$ is inert in $\mathbf{k}$, we require that

$$(5.20) \qquad (\Delta, -pp_0)_v = \begin{cases} -1 & \text{if } v = p \text{ or } \infty, \text{ and} \\ 1 & \text{otherwise.} \end{cases}$$

If $p$ is ramified in $\mathbf{k}$, we require that

$$(5.21) \qquad (\Delta, -p_0)_v = \begin{cases} -1 & \text{if } v = p \text{ or } \infty, \text{ and} \\ 1 & \text{otherwise.} \end{cases}$$

Let

$$(5.22) \qquad \kappa_p = \begin{cases} pp_0 & \text{if } p \text{ is inert in } \mathbf{k}, \text{ and} \\ p_0 & \text{if } p \text{ is ramified in } \mathbf{k}. \end{cases}$$

Then $\mathbb{B}_p$ is the cyclic algebra $(\Delta, -\kappa_p)$ and thus can be written in the form

$$(5.23) \qquad \mathbb{B}_p = \mathbf{k} \oplus \mathbf{k} \cdot y$$

for an element $y$ with $\operatorname{tr}(y) = 0$, $y^2 = -\kappa_p$, and $a \cdot y = y \cdot \bar{a}$ for all $a \in \mathbf{k}$.

Let $R = \mathcal{O}_{\mathbf{k}} \oplus \mathcal{O}_{\mathbf{k}} \cdot y$, so that $R$ is an order in $\mathbb{B}_p$ with $R \cap \mathbf{k} = \mathcal{O}_{\mathbf{k}}$. Suppose that $\mathcal{O}$ is a maximal order of $\mathbb{B}_p$ containing $R$. If $x_1 \in R$ and $x_2 \in \mathcal{O}$, then $\operatorname{tr}(x_1 x_2) \in \mathbb{Z}$. Using the coordinates from (5.23), this gives

$$(5.24) \qquad \operatorname{tr}((\alpha_1 + \beta_2 y)(\alpha_2 + \beta_2 y)) = \operatorname{tr}_{\mathbf{k}}(\alpha_1 \alpha_2 - \kappa_p \beta_1 \bar{\beta}_2) \in \mathbb{Z}.$$

Therefore

$$(5.25) \qquad R = \mathcal{O}_{\mathbf{k}} \oplus \mathcal{O}_{\mathbf{k}} \cdot y \subset \mathcal{O} \subset \mathcal{D}^{-1} \oplus \kappa_p^{-1} \mathcal{D}^{-1} \cdot y,$$

where $\mathcal{D}^{-1} = \sqrt{\Delta}^{-1} \mathcal{O}_{\mathbf{k}}$ is the inverse different of $\mathbf{k}$. Setting

$$(5.26) \qquad V(\mathcal{O}) = V^{(p)} \cap \mathcal{O},$$

we conclude that

$$\mathcal{O}_{\mathbf{k}} \cdot y \subset V(\mathcal{O}) \subset \kappa_p^{-1} \mathcal{D}^{-1} \cdot y. \tag{5.27}$$

The conditions imposed on $p_0$ imply that $p_0$ splits. Let $p_0 \mathcal{O}_{\mathbf{k}} = \wp_0 \bar{\wp}_0$ and let $p\mathcal{O}_{\mathbf{k}} = \wp^{e_p}$. Then, it is easily checked that

$$R' := \mathcal{O}_{\mathbf{k}} \oplus \wp_0^{-1} y \quad \text{and} \quad R'' := \mathcal{O}_{\mathbf{k}} \oplus \bar{\wp}_0^{-1} y \tag{5.28}$$

are orders in $\mathbb{B}_p$ containing $R$. Note that $R'$ and $R''$ cannot be contained in the same order $\mathcal{O}$. Indeed, if $R' + R'' \subset \mathcal{O}$, then $(\wp_0^{-1} + \bar{\wp}_0^{-1}) \cdot y = p_0^{-1} \mathcal{O}_{\mathbf{k}} \cdot y \subset \mathcal{O}$. But $Q(p_0^{-1} y) = -(p_0^{-1} y)^2 = p_0^{-2} \kappa_p \notin \mathbb{Z}$.

Suppose that $\mathcal{O}$ is a maximal order containing $\mathcal{O}_{\mathbf{k}} \oplus \wp_0^{-1} y$. If $x = \kappa_p^{-1} \sqrt{\Delta}^{-1} \beta y \in V(\mathcal{O})$, with $\beta \in \mathcal{O}_{\mathbf{k}}$, then

$$Q(x) = -x^2 = \kappa_p^{-1} \Delta^{-1} N(\beta) \in \mathbb{Z}. \tag{5.29}$$

This forces

$$\beta \in \wp_0^{-1} \kappa_p \sqrt{\Delta} \mathcal{O}_{\mathbf{k}}, \quad \text{or} \quad \beta \in \bar{\wp}_0^{-1} \kappa_p \sqrt{\Delta} \mathcal{O}_{\mathbf{k}}, \tag{5.30}$$

so that

$$x \in \wp_0^{-1} y \quad \text{or} \quad x \in \bar{\wp}_0^{-1} y. \tag{5.31}$$

Therefore,

$$V(\mathcal{O}) = \wp_0^{-1} \cdot y. \tag{5.32}$$

Fix a maximal order $\mathcal{O}$ containing $R'$, and note that $\mathcal{O} \cap \mathbf{k} = \mathcal{O}_{\mathbf{k}}$, i.e., that $\mathcal{O}_{\mathbf{k}}$ is optimally embedded in $\mathcal{O}$, in Eichler's terminology, [**8**]. By the Chevalley-Hasse-Noether Theorem (comp. [**8**], Satz 7), it follows that if $\mathcal{O}'$ is any maximal order in $\mathbb{B}_p$ in which $\mathcal{O}_{\mathbf{k}}$ is optimally embedded, then there exists a finite idèle $g \in \mathbf{k}_{\mathbb{A}_f}^\times$, such that

$$\mathcal{O}' = g(\mathcal{O} \otimes_{\mathbb{Z}} \hat{\mathbb{Z}}) g^{-1} \cap \mathbb{B}_p. \tag{5.33}$$

If

$$\mathfrak{a} = \left(g\mathcal{O}_{\mathbf{k}} \otimes_{\mathbb{Z}} \hat{\mathbb{Z}}\right) \cap \mathbf{k}, \tag{5.34}$$

then

$$V(\mathcal{O}') = \mathfrak{a}\bar{\mathfrak{a}}^{-1} \wp_0^{-1} y. \tag{5.35}$$

Conversely, any ideal $\mathfrak{a}$ arises in this way.

Now suppose that $\mathcal{O}' = \text{End}(E)$ for a supersingular elliptic curve $(E, \iota) \in \mathcal{M}(\overline{\kappa_p(\wp)})$ over $\bar{\mathbb{F}}_p$ with complex multiplication by $\mathcal{O}_{\mathbf{k}}$. Note that $\mathcal{O}_{\mathbf{k}}$ is optimally embedded in $\mathcal{O}'$ and that

$$V(E, \iota) = \mathcal{O}' \cap \mathbf{k}y = V(\mathcal{O}'). \tag{5.36}$$

**Proposition 5.13.** *Fix an auxillary prime $p_0$ as above. Then there is an element $y_0$ in $V$ with $Q(y_0) = \kappa_p$, and the lattices $V(E_i, \iota_i)$ in Theorem 5.12 have the form*
$$V(E_i, \iota_i) \simeq \mathfrak{a}_i \bar{\mathfrak{a}}_i^{-1} \wp_0^{-1} \cdot y_0,$$
*as the $\mathfrak{a}_i$ run over representatives for the ideal classes of $\mathbf{k}$. In particular, the isomorphism class of $V(E_i, \iota_i)$ runs over the genus $[[\wp_0]]$ of the ideal $\wp_0$.*

**Remark.** The genus $[[\wp_0]]$ is independent of the choice of $\wp_0$ and is characterized by the values of the genus characters. Recall that a basis for the characters of the group of genera are given as follows, [**12**]. Let $q_i$, $1 \leq i \leq t$ be the primes dividing $\Delta$. Let

(5.37)
$$\chi_i(\mathfrak{a}) = (\Delta, N(\mathfrak{a}))_{q_i},$$

where $( \, , \, )_{q_i}$ is the Hilbert symbol for $\mathbb{Q}_{q_i}$. Note that, if $x \in \mathbf{k}^\times$, then $(\Delta, N(x))_{q_i} = 1$, and that

(5.38)
$$(\Delta, N(\mathfrak{a}) N(\mathfrak{b})^2)_{q_i} = (\Delta, N(\mathfrak{a}))_{q_i}.$$

Thus $\chi_i$ depends only on the ideal class of $\mathfrak{b}$ modulo squares, i.e., defines a character of the group of genera $Cl(\mathbf{k})/2Cl(\mathbf{k})$. Any $t-1$ of the $\chi_i$'s give a basis for the characters of this group. The conditions (5.20) and (5.21) on the auxillary prime $p_0$ imply that

(5.39) $\quad \chi_i(\wp_0) = (\Delta, p_0)_{q_i} = \begin{cases} (\Delta, -p)_{q_i} & \text{if } p \text{ is inert,} \\ (\Delta, -1)_{q_i} & \text{if } p \text{ is ramified and } q_i \neq p, \text{ and} \\ -(\Delta, -1)_{q_i} & \text{if } p \text{ is ramified and } q_i = p. \end{cases}$

The genus of $\wp_0$ is determined by these conditions.

**Corollary 5.14.** *(i) The rational quadratic space $V^{(p)}$ is given by*
$$(V^{(p)}, Q) \simeq (\mathbf{k}, \kappa_p N),$$
*where $\kappa_p$ is given by (5.22).*
*(ii) Let $C_i$ be the ideal class*
$$C_i = [\mathfrak{a}_i^{-1} \bar{\mathfrak{a}}_i \, \wp_0] = [\mathfrak{a}_i]^2 [\wp_0].$$
*Then*
$$\left| \{ x \in V(E_i, \iota_i) \mid Q(x) = t \} \right| = |\mathcal{O}_\mathbf{k}^\times| \cdot \left| \{ \mathfrak{c} \subset \mathcal{O}_\mathbf{k} \mid [\mathfrak{c}] = C_i, \text{ and } N(\mathfrak{c}) = p_0 t / \kappa_p \} \right|.$$

Note that

$$(5.40) \qquad p_0/\kappa_p = \begin{cases} p & \text{if } p \text{ is inert, and} \\ 1 & \text{if } p \text{ is ramified.} \end{cases}$$

*Proof.* We just note that, if $x = \beta \cdot y \in V(E_i, \iota_i)$, with $\beta \in \mathfrak{a}_i \bar{\mathfrak{a}}_i^{-1} \wp_0^{-1}$, then $t = Q(x) = \kappa_p N(\beta)$. The ideal $\mathfrak{c} = \beta \mathfrak{a}_i^{-1} \bar{\mathfrak{a}}_i \wp_0$ is integral, lies in the ideal class $C_i$, and has

$$(5.41) \qquad N(\mathfrak{c}) = N(\beta) p_0 = t p_0/\kappa_p.$$

□

Note that, by genus theory, as $i$ runs from $1$ to $h_{\boldsymbol{k}}$, $[\mathfrak{a}_i]$ runs over the ideal class group and $C_i$ runs over $(2Cl(\boldsymbol{k}))[\wp_0] = [[\wp_0]]$, the genus of $\wp_0$, with each class in $[[\wp_0]]$ occurring $2^{t-1}$ times. Here again $t$ is the number of rational primes which ramify in $\boldsymbol{k}$.

Finally, we specialize to the case of prime discriminant considered in earlier sections of this paper, i.e., let $\boldsymbol{k} = \mathbb{Q}(\sqrt{-q})$ with $q > 3$ a prime congruent to $3$ modulo $4$. In this case, the class number of $\boldsymbol{k}$ is odd, there is only one genus, and the $C_i$'s run over all ideal classes. Combining the results above, we obtain:

**Theorem 5.15.**

$$\deg(Z(t)) = 2 \log(q) \cdot (\mathrm{ord}_q(t) + 1) \cdot \rho(t) + 2 \sum_p \log(p) \cdot (\mathrm{ord}_p(t) + 1) \cdot \rho(t/p).$$

*Here $p$ runs over the primes which are inert in $\boldsymbol{k}$.*

This is precisely the expression given in Theorem 3 of the introduction.

### §6. Variants.

In this section, we will sketch a few variations and extensions of the results of the previous sections. However, none of them touch on the most tantalizing problem of establishing a *direct connection* between the Fourier coefficients and the degrees of special cycles.

The central derivatives of *all* of the incoherent Eisenstein series constructed in section 1 should have an arithmetic interpretation analogous to Theorem 3. A full description of such a result would require several extensions of what we have done above, which is specialized in several respects.

First of all, we have restricted to the case of a prime discriminant $\Delta = -q$. This eliminates complications involving genus theory. Of course, the moduli problem of section 5 has been set up to allow elliptic curves with complex multiplication by $\mathcal{O}_{\boldsymbol{k}}$ where $\boldsymbol{k} = \mathbb{Q}(\sqrt{\Delta})$ for an arbitrary fundamental discriminant $\Delta < 0$. On the analytic side, the machinery of section 1 and 2 provides corresponding incoherent Eisenstein series, whose Fourier expansions can be computed by the methods of section 2. It should be then a routine matter to work out the analogues of Theorems 1, 2 and 3 in this case, except that the group of genera will now play a nontrivial role.

Second, in section 2, we have considered only those Eisenstein series built from sections which are defined by the characteristic functions of the completions of $\mathcal{O}_{\boldsymbol{k}}$. More general sections could be considered at the cost of (i) more elaborate calculations of Whittaker functions and their derivatives used in computing the Fourier expansion of the central derivative, and (ii) incorporation of a level structure in the moduli problem.

Third, in section 2, we have restricted to the simplest type of incoherent collection $\mathcal{C}$, i.e., a collection $\mathcal{C}$ which is obtained by considering the global quadratic space $(V,Q) = (\boldsymbol{k}, -N)$ and switching the signature from $(0,2)$ to $(2,0)$. Thus, at each finite place, $\mathcal{C}_p = (\boldsymbol{k}_p, -N)$. The most general collection $\mathcal{C}$ is obtained by switching the signature of the space $(V,Q) = (\boldsymbol{k}, \kappa N)$ where $\kappa \in \mathbb{Q}^\times$. In this situation, a slightly more elaborate moduli problem is necessary. We give a brief sketch of this.

Let $B$ be a indefinite quaternion algebra over $\mathbb{Q}$, with a fixed embedding $i: \boldsymbol{k} \hookrightarrow B$. Fix a maximal order $\mathcal{O}_B$, in $B$, a positive involution $*$ of $B$ stabilizing $\mathcal{O}_B$, and assume that $i(\mathcal{O}_{\boldsymbol{k}}) \subset \mathcal{O}_B$.

We consider polarized abelian surfaces $A$ with an action $\iota : \mathcal{O}_B \otimes \mathcal{O}_{\boldsymbol{k}} \to \mathrm{End}(A)$ compatible with the positive involution $b \otimes a \mapsto (b \otimes a)^* = b^* \otimes \bar{a}$. For such an $(A, \iota)$, let

(6.1) $$V(A, \iota) = \{ y \in \mathrm{End}(A) \mid y\iota(b \otimes a) = \iota(b \otimes \bar{a})y \}$$

be the space of special endomorphisms. Observe that

(6.2) $$M_2(\boldsymbol{k}) \simeq B \otimes \boldsymbol{k} \hookrightarrow \mathrm{End}^0(A),$$

so that $A$ is isogenous to a product $A_0^2$ for an elliptic curve $A_0$ with complex multiplication. Thus

(6.3) $$\mathrm{End}^0(A) \simeq M_2(\mathrm{End}^0(A_0)).$$

If $A_0$ is *not* a supersingular elliptic curve in characteristic $p$, for $p$ not split in $\pmb{k}$, then $\mathrm{End}^0(A) \simeq M_2(\pmb{k})$ and $V^0(A, \iota) = V(A, \iota) \otimes_{\mathbb{Z}} \mathbb{Q} = 0$. If $A_0$ is supersingular in characteristic $p$, then

$$\text{(6.4)} \qquad \mathrm{End}^0(A) \simeq M_2(\mathbb{B}_p),$$

where $\mathbb{B}_p = \mathrm{End}^0(A_0)$ is, as before, the definite quaternion algebra ramified at infinity and $p$. The commutator of $\iota(B \otimes 1)$ in $\mathrm{End}^0(A)$ is then isomorphic to $B'$, the definite quaternion algebra over $\mathbb{Q}$ with

$$\text{(6.5)} \qquad \mathrm{inv}_\ell(B') = \begin{cases} \mathrm{inv}_\ell(B) & \text{if } \ell \neq p, \text{ and} \\ -\mathrm{inv}_\ell(B') & \text{if } \ell = p. \end{cases}$$

Note that $\iota$ yields an embedding of $\pmb{k}$ into $B'$, and we have a decomposition, analogous to (5.16),

$$\text{(6.6)} \qquad B' = \pmb{k} \oplus V^0(A, \iota),$$

where $V^0(A, \iota) = \pmb{k} \cdot y$ for an nonzero element $y \in B'$ with $ya = \bar{a}y$ for all $a \in \pmb{k}$. The quadratic form on $V^0(A, \iota)$ is then

$$\text{(6.7)} \qquad Q(\beta y) = -(\beta y)^2 = -N(\beta) y^2,$$

so that the binary quadratic space $(V^0(A, \iota), Q)$ is isomorphic to $(\pmb{k}, \kappa N)$ where $\kappa = -y^2$. Note that $B'$ is then isomorphic to the cyclic algebra $(\Delta, -\kappa)$. An analysis of the corresponding integral moduli problem like that of section 5 should be possible and give the desired generalization of the results above. The case $B = M_2(\mathbb{Q})$ will reduce precisely to the situation of section 5.

Returning to Theorem 3 in the case $\pmb{k} = \mathbb{Q}(\sqrt{-q})$, from the point of view of Arakelov theory, it is possible to given an interpretation of the negative Fourier coefficients of $\phi$ as degrees of cycles.

First recall, [**22**], that an Arakelov divisor for $M = \mathrm{Spec}(\mathcal{O}_H)$ is an expression of the form

$$\text{(6.8)} \qquad D = \sum_\lambda^{h_{\pmb{k}}} r_\lambda \lambda + \sum_{\mathfrak{p}} n_{\mathfrak{p}} \mathfrak{p},$$

where $\mathfrak{p}$ runs over the nonzero prime ideals of $\mathcal{O}_H$, $\lambda$ runs over the complex places of $H$, the $n_{\mathfrak{p}}$'s are integers (almost all zero), and the $r_\lambda$'s are real numbers. Let

35$\mathrm{Div}_c(M)$ be the group of Arakelov divisors, and let $\mathrm{Pic}_c(M)$ be its quotient by the group of principal Arakelov divisors, [**22**].

We can view the cycles $Z(t)$ defined in section 2 as elements of $\mathrm{Div}_c(M)$ as follows. Fix the embedding of $\boldsymbol{k}$ into $\mathbb{C}$, where $\sqrt{-q}$ has positive imaginary part, and let $\omega = (1+\sqrt{-q})/2$ so that $\mathcal{O}_{\boldsymbol{k}} = \mathbb{Z}+\mathbb{Z}\omega$. Let

$$(6.9) \qquad j_0 = j(\omega) = j(\mathcal{O}_{\boldsymbol{k}}) \in \mathcal{O}_H$$

be the corresponding singular value of the $j$-function.

Let $p$ be a rational prime which does not split in $\boldsymbol{k}$ and let $p\mathcal{O}_{\boldsymbol{k}} = \wp^{e_p}$. Then $\wp$ is principal and hence splits completely in $H$. For any prime $\mathfrak{p}$ of $\mathcal{O}_H$ above $\wp$, the image of $j_0$ under the reduction map $\mathcal{O}_H \to \mathcal{O}_H/\mathfrak{p} = \mathbb{F}_{p^{3-e_p}}$ is the $j$-invariant of an elliptic curve $E_{\mathfrak{p}}$ over $\bar{\mathbb{F}}_p$, unique up to isomorphism, with complex multiplication by $\mathcal{O}_{\boldsymbol{k}}$. For each such $p$, this gives a bijective map from the primes over $\wp$ to the fiber $M(\bar{\mathbb{F}}_p)$. Then, for $t > 0$,

$$(6.10) \qquad Z(t) = \sum_p \sum_{\mathfrak{p}\mid\wp} n_{\mathfrak{p}}(t)\mathfrak{p},$$

where

$$(6.11) \qquad n_{\mathfrak{p}}(t) = \nu_p(t) \cdot \big|\{x \in V(E_{\mathfrak{p}},\iota) \mid Q(x) = t\}\big|.$$

Here $V(E_{\mathfrak{p}},\iota)$ is the lattice of special endomorphisms of $E_{\mathfrak{p}}$, as in Definition 5.7, and $\nu_p(t)$ is the integer given in Theorem 5.11. Let

$$(6.11) \qquad C_{\mathfrak{p}} = [V(E_{\mathfrak{p}},\iota)]^{-1},$$

where $[V(E_{\mathfrak{p}},\iota)] \in \mathrm{Pic}(\mathcal{O}_{\boldsymbol{k}}) \simeq Cl(\boldsymbol{k})$ is the class of the rank one $\mathcal{O}_{\boldsymbol{k}}$-module $[V(E_{\mathfrak{p}},\iota)]$. Using Corollary 5.14, (6.10) becomes

$$(6.12) \qquad n_{\mathfrak{p}} = \nu_p(t) \cdot \big|\{\mathfrak{c} \subset \mathcal{O}_{\boldsymbol{k}} \mid [\mathfrak{c}] = C_{\mathfrak{p}} \text{ and } N(\mathfrak{c}) = t/\kappa_p\}\big|,$$

where

$$(6.13) \qquad \kappa_p = \begin{cases} p & \text{if } p \text{ is inert in } \boldsymbol{k}, \text{ and} \\ 1 & \text{if } p = q. \end{cases}$$

We view $Z(t)$ as an element of $\mathrm{Div}_c(M)$ with zero archimedean component.

To define an element $Z(t) \in \mathrm{Div}_c(M)$ when $t < 0$, we use a construction based on an idea explained to us by Dick Gross. For a point $(E,\iota)$ of $\mathcal{M}(\mathbb{C})$, let $E^{\mathrm{top}}$ be the underlying real torus. By analogy with Definition 5.7, we define the space of special endomorphisms

$$V(E, \iota) = \{y \in \mathrm{End}(E^{\mathrm{top}}) \mid y\iota(a) = \iota(\bar{a})y\}. \tag{6.14}$$

If we write $E = \mathbb{C}/L$ for a lattice $L$ with a chosen basis, then

$$\mathrm{End}(E^{\mathrm{top}}) \simeq M_2(\mathbb{Z}), \tag{6.15}$$

and so, by analogy with the case of a supersingular elliptic curve in characteristic $p$, $\mathrm{End}(E^{\mathrm{top}})$ is isomorphic to a maximal order in the quaternion algebra $M_2(\mathbb{Q})$. It is not difficult to check that $V(E, \iota)$ is a rank one $\mathcal{O}_{\boldsymbol{k}}$-module with a $\mathbb{Z}$-valued quadratic form given by $Q(y) = -y^2$.

For each embedding $\lambda$ of $H = \boldsymbol{k}(j)$ into $\mathbb{C}$, let $E_\lambda$ be the elliptic curve over $\mathbb{C}$ with $j$-invariant $\lambda(j_0)$. This curve has complex multiplication by $\mathcal{O}_{\boldsymbol{k}}$ and space of special endomorphism, in the sense just defined, $V(E_\lambda, \iota)$. Let $C_\lambda = [V(E_\lambda, \iota)]^{-1} \in \mathrm{Pic}(\mathcal{O}_{\boldsymbol{k}})$. Then, for $t \in \mathbb{Z}_{<0}$ and for $\tau = u + iv$ in the upper half plane, we associate the Arakelov divisor

$$Z(t) = Z(t, v) := \sum_\lambda r_\lambda \lambda, \tag{6.16}$$

where

$$\begin{aligned} r_\lambda(t) = r_\lambda(t, v) &= 2\beta_1(4\pi|t|v) \cdot |\{y \in V(E_\lambda, \iota) \mid Q(y) = t\} \\ &= 2\beta_1(4\pi|t|v) \cdot |\{\mathfrak{c} \subset \mathcal{O}_{\boldsymbol{k}} \mid [\mathfrak{c}] = C_\lambda \text{ and } N(\mathfrak{c}) = -t\}|. \end{aligned} \tag{6.17}$$

Here $\beta_1$ is the exponential integral (0.6). Thus we have defined for all $t \neq 0$ an Arakelov divisor with

$$\deg(Z(t)) = \begin{cases} a_t(\phi) & \text{if } t > 0, \text{ and} \\ a_t(\phi, v) & \text{if } t < 0, \end{cases} \tag{6.18}$$

and hence

$$\phi(\tau) = a_0(\phi, v) + \sum_{t \neq 0} \deg(Z(t))\, q^t, \tag{6.19}$$

extending (0.8). It remains to find a natural definition of an Arakelov divisor $Z(0) = Z(0, v)$ such that

$$\deg(Z(0)) = a_0(\phi, v) = -h_{\boldsymbol{k}}\left(\log(q) + \log(v) + 2\frac{\Lambda'(1, \chi)}{\Lambda(1, \chi)}\right). \tag{6.20}$$

Supposing this done, one would be tempted to consider the generating series

$$\sum_{t \in \mathbb{Z}} [Z(t)]\, q^t \tag{6.21}$$

where $[Z(t)]$ is the image of $Z(t)$ in $\mathrm{Pic}_c(M) \otimes_{\mathbb{Z}} \mathbb{C}$. We do not know if such a series has a meaning (e.g., converges).

S. Kudla  
Department of Mathematics  
University of Maryland  
College Park MD 20742  
USA

M. Rapoport  
Mathematisches Institut  
der Universität zu Köln  
Weyertal 86-90  
50931 Köln  
Germany

T. Yang  
Department of Mathematics  
SUNY at Stony Brook  
Stony Brook, NY 11794  
USA